\documentclass[12pt]{article}
\usepackage{full page, 
amssymb, amscd, graphicx, 
}
    \title{{\bf Logarithmic
tensor category theory, VII: Convergence and extension
properties and applications to expansion for intertwining
maps}}
    \author{Yi-Zhi Huang, James Lepowsky and Lin Zhang}
    \date{}

\begin{document}
    \bibliographystyle{alpha}
    \maketitle

    \newtheorem{rema}{Remark}[section]
    \newtheorem{propo}[rema]{Proposition}
    \newtheorem{theo}[rema]{Theorem}
   \newtheorem{defi}[rema]{Definition}
    \newtheorem{lemma}[rema]{Lemma}
    \newtheorem{corol}[rema]{Corollary}
     \newtheorem{exam}[rema]{Example}
\newtheorem{assum}[rema]{Assumption}
     \newtheorem{nota}[rema]{Notation}
        \newcommand{\ba}{\begin{array}}
        \newcommand{\ea}{\end{array}}
        \newcommand{\be}{\begin{equation}}
        \newcommand{\ee}{\end{equation}}
        \newcommand{\bea}{\begin{eqnarray}}
        \newcommand{\eea}{\end{eqnarray}}
        \newcommand{\nno}{\nonumber}
        \newcommand{\nn}{\nonumber\\}
        \newcommand{\lbar}{\bigg\vert}
        \newcommand{\p}{\partial}
        \newcommand{\dps}{\displaystyle}
        \newcommand{\bra}{\langle}
        \newcommand{\ket}{\rangle}
 \newcommand{\res}{\mbox{\rm Res}}
\newcommand{\wt}{\mbox{\rm wt}\;}
\newcommand{\swt}{\mbox{\scriptsize\rm wt}\;}
 \newcommand{\pf}{{\it Proof}\hspace{2ex}}
 \newcommand{\epf}{\hspace{2em}$\square$}
 \newcommand{\epfv}{\hspace{1em}$\square$\vspace{1em}}
        \newcommand{\ob}{{\rm ob}\,}
        \renewcommand{\hom}{{\rm Hom}}
\newcommand{\C}{\mathbb{C}}
\newcommand{\R}{\mathbb{R}}
\newcommand{\Z}{\mathbb{Z}}
\newcommand{\N}{\mathbb{N}}
\newcommand{\A}{\mathcal{A}}
\newcommand{\Y}{\mathcal{Y}}
\newcommand{\Arg}{\mbox{\rm Arg}\;}
\newcommand{\comp}{\mathrm{COMP}}
\newcommand{\lgr}{\mathrm{LGR}}

\newcommand{\dlt}[3]{#1 ^{-1}\delta \bigg( \frac{#2 #3 }{#1 }\bigg) }

\newcommand{\dlti}[3]{#1 \delta \bigg( \frac{#2 #3 }{#1 ^{-1}}\bigg) }

 \makeatletter
\newlength{\@pxlwd} \newlength{\@rulewd} \newlength{\@pxlht}
\catcode`.=\active \catcode`B=\active \catcode`:=\active \catcode`|=\active
\def\sprite#1(#2,#3)[#4,#5]{
   \edef\@sprbox{\expandafter\@cdr\string#1\@nil @box}
   \expandafter\newsavebox\csname\@sprbox\endcsname
   \edef#1{\expandafter\usebox\csname\@sprbox\endcsname}
   \expandafter\setbox\csname\@sprbox\endcsname =\hbox\bgroup
   \vbox\bgroup
  \catcode`.=\active\catcode`B=\active\catcode`:=\active\catcode`|=\active
      \@pxlwd=#4 \divide\@pxlwd by #3 \@rulewd=\@pxlwd
      \@pxlht=#5 \divide\@pxlht by #2
      \def .{\hskip \@pxlwd \ignorespaces}
      \def B{\@ifnextchar B{\advance\@rulewd by \@pxlwd}{\vrule
         height \@pxlht width \@rulewd depth 0 pt \@rulewd=\@pxlwd}}
      \def :{\hbox\bgroup\vrule height \@pxlht width 0pt depth
0pt\ignorespaces}
      \def |{\vrule height \@pxlht width 0pt depth 0pt\egroup
         \prevdepth= -1000 pt}
   }
\def\endsprite{\egroup\egroup}
\catcode`.=12 \catcode`B=11 \catcode`:=12 \catcode`|=12\relax
\makeatother

\def\hboxtr{\FormOfHboxtr} 
\sprite{\FormOfHboxtr}(25,25)[0.5 em, 1.2 ex] 

:BBBBBBBBBBBBBBBBBBBBBBBBB |
:BB......................B |
:B.B.....................B |
:B..B....................B |
:B...B...................B |
:B....B..................B |
:B.....B.................B |
:B......B................B |
:B.......B...............B |
:B........B..............B |
:B.........B.............B |
:B..........B............B |
:B...........B...........B |
:B............B..........B |
:B.............B.........B |
:B..............B........B |
:B...............B.......B |
:B................B......B |
:B.................B.....B |
:B..................B....B |
:B...................B...B |
:B....................B..B |
:B.....................B.B |
:B......................BB |
:BBBBBBBBBBBBBBBBBBBBBBBBB |

\endsprite

\def\shboxtr{\FormOfShboxtr} 
\sprite{\FormOfShboxtr}(25,25)[0.3 em, 0.72 ex] 

:BBBBBBBBBBBBBBBBBBBBBBBBB |
:BB......................B |
:B.B.....................B |
:B..B....................B |
:B...B...................B |
:B....B..................B |
:B.....B.................B |
:B......B................B |
:B.......B...............B |
:B........B..............B |
:B.........B.............B |
:B..........B............B |
:B...........B...........B |
:B............B..........B |
:B.............B.........B |
:B..............B........B |
:B...............B.......B |
:B................B......B |
:B.................B.....B |
:B..................B....B |
:B...................B...B |
:B....................B..B |
:B.....................B.B |
:B......................BB |
:BBBBBBBBBBBBBBBBBBBBBBBBB |

\endsprite


\begin{abstract}
This is the seventh part in a series of papers in which we introduce
and develop a natural, general tensor category theory for suitable
module categories for a vertex (operator) algebra.  In this paper
(Part VII), we give sufficient conditions for the existence of the
associativity isomorphisms.
\end{abstract}


\tableofcontents
\vspace{2em}

In this paper, Part VII of a series of eight papers on logarithmic
tensor category theory, we give sufficient conditions for the
existence of the associativity isomorphisms.  The sections, equations,
theorems and so on are numbered globally in the series of papers
rather than within each paper, so that for example equation (a.b) is
the b-th labeled equation in Section a, which is contained in the
paper indicated as follows: In Part I \cite{HLZ1}, which contains
Sections 1 and 2, we give a detailed overview of our theory, state our
main results and introduce the basic objects that we shall study in
this work.  We include a brief discussion of some of the recent
applications of this theory, and also a discussion of some recent
literature.  In Part II \cite{HLZ2}, which contains Section 3, we
develop logarithmic formal calculus and study logarithmic intertwining
operators.  In Part III \cite{HLZ3}, which contains Section 4, we
introduce and study intertwining maps and tensor product bifunctors.
In Part IV \cite{HLZ4}, which contains Sections 5 and 6, we give
constructions of the $P(z)$- and $Q(z)$-tensor product bifunctors
using what we call ``compatibility conditions'' and certain other
conditions.  In Part V \cite{HLZ5}, which contains Sections 7 and 8,
we study products and iterates of intertwining maps and of logarithmic
intertwining operators and we begin the development of our analytic
approach.  In Part VI \cite{HLZ6}, which contains Sections 9 and 10,
we construct the appropriate natural associativity isomorphisms
between triple tensor product functors.  The present paper, Part VII,
contains Section 11.  In Part VIII \cite{HLZ8}, which contains Section
12, we construct braided tensor category structure.

\paragraph{Acknowledgments}
The authors gratefully
acknowledge partial support {}from NSF grants DMS-0070800 and
DMS-0401302.  Y.-Z.~H. is also grateful for partial support {}from NSF
grant PHY-0901237 and for the hospitality of Institut des Hautes 
\'{E}tudes Scientifiques in the fall of 2007.

\renewcommand{\theequation}{\thesection.\arabic{equation}}
\renewcommand{\therema}{\thesection.\arabic{rema}}
\setcounter{section}{10}
\setcounter{equation}{0}
\setcounter{rema}{0}

\section{The convergence and extension properties 
and differential equations}

In the construction of the associativity isomorphisms we have needed,
and assumed, the convergence and expansion conditions for intertwining
maps in ${\cal C}$.  In this section we will follow \cite{tensor4} to
give certain sufficient conditions for a category ${\cal C}$ to have
these properties. In Section 11.1, by exhibiting explicitly the shapes
of products and iterates of logarithmic intertwining operators as
analytic functions, we give what we call, as in \cite{tensor4}, the
``convergence and extension properties'' and show that they imply the
convergence and expansion conditions for intertwining maps in ${\cal
C}$. In our general setting, these properties in general involve both
the logarithm function and the abelian group gradings.  These
properties are formulated {\it globally} using analytic functions
defined on the regions $|z_{1}|>|z_{2}|>0$ and
$|z_{2}|>|z_{1}-z_{2}|>0$, while the expansion condition is formulated
{\it locally} near a point $(z_1,z_2)$ in the intersection of these
regions, that is, such that $|z_{1}|>|z_{2}|>|z_{1}-z_{2}|>0$.  In
Section 11.2, we show that the proofs in \cite{diff-eqn}, deriving and
using differential equations, can be adapted to generalize the results
of \cite{diff-eqn} to results in the logarithmic generality. In
particular, we will see that two purely algebraic conditions, the
``$C_{1}$-cofiniteness condition'' and the
``quasi-finite-dimensionality condition,'' for all objects of
$\mathcal{C}$ imply the convergence and extension properties and thus
also imply the convergence and expansion conditions for intertwining
maps in ${\cal C}$.

\subsection{The convergence and extension properties}

Given objects $W_1$, $W_2$, $W_3$, $W_4$, $M_1$ and $M_2$ of the
category ${\cal C}$, let ${\cal Y}_1$, ${\cal Y}_2$, ${\cal Y}^1$ and
${\cal Y}^2$ be logarithmic intertwining operators of types
${W_4}\choose {W_1M_1}$, ${M_1}\choose {W_2W_3}$, ${W_4}\choose
{M_2W_3}$ and ${M_2}\choose {W_1W_2}$, respectively.  We now consider
certain natural conditions on the product of ${\cal Y}_1$ and ${\cal
Y}_2$ and on the iterate of ${\cal Y}^1$ and ${\cal Y}^2$.  These
conditions require that the product of ${\cal Y}_1$ and ${\cal Y}_2$
(respectively, the iterate of ${\cal Y}^1$ and ${\cal Y}^2$) be
absolutely convergent in the region $|z_{1}|>|z_{2}|>0$ (respectively,
$|z_{2}|>|z_{1}-z_{2}|>0$) and that it can be analytically extended to
a function of the same shape as (finite) sums of iterates
(respectively, products) of logarithmic intertwining operators in the
region $|z_{2}|>|z_{1}-z_{2}|>0$ (respectively, $|z_{1}|>|z_{2}|>0$)
with finitely-generated ``lower bounded doubly-graded generalized $V$-modules'' as
intermediate generalized $V$-modules (see Definition \ref{d-g-lower-b}
below).  In fact, when the intermediate generalized $V$-module $M_{2}$
is a finitely-generated lower bounded doubly-graded generalized $V$-module, the sum
defining the iterate (\ref{c-e-p-6}) in the region
$|z_{2}|>|z_{1}-z_{2}|>0$ is a branch of a multivalued analytic
function of the form (\ref{c-e-p-3}), and analogously, when the
intermediate generalized $V$-module $M_{1}$ is a finitely-generated lower bounded
doubly-graded generalized $V$-module, the sum defining the product
(\ref{c-e-p-2}) in the region $|z_{1}|>|z_{2}|>0$ is a branch of a
multivalued analytic function of the form (\ref{c-e-p-7}). In the
following conditions, we recall from the beginning of Section 7, in
particular, (\ref{convp}), (\ref{convi}), (\ref{4itm}) and
(\ref{4prm}), the meaning of the absolute convergence of
(\ref{c-e-p-2}) and (\ref{c-e-p-6}) below; as in (\ref{iterabbr}) and
(\ref{prodabbr}), we are taking $p=0$ in the notation of
(\ref{im:f(z)}).

\begin{description}

\item[Convergence and extension property for products] For any
$\beta\in \tilde{A}$, there exists an integer $N_{\beta}$ depending
only on ${\cal Y}_1$, ${\cal Y}_2$ and $\beta$, and for any
weight-homogeneous elements $w_{(1)}\in (W_1)^{(\beta_{1})}$ and
$w_{(2)}\in (W_2)^{(\beta_{2})}$ ($\beta_{1}, \beta_{2} \in
\tilde{A}$) and any $w_{(3)}\in W_3$ and $w'_{(4)}\in W'_4$ such that
\begin{equation}\label{c-e-p-0}
\beta_{1}+\beta_{2}=-\beta, 
\end{equation}
there exist $M\in{\mathbb N}$, $r_{k}, s_{k}\in {\mathbb R}$, $i_{k},
j_{k}\in \N$, $k=1,\dots,M$; $K \in \Z_{+}$ independent of $w_{(1)}$ and
$w_{(2)}$ such that each $i_{k} < K$; and analytic functions $f_{k}(z)$ on
$|z|<1$, $k=1, \dots, M$, satisfying
\begin{equation}\label{c-e-p-1}
\wt w_{(1)}+\wt w_{(2)}+s_{k}>N_{\beta}, \;\;\;k=1, \dots, M,
\end{equation}
such that
\begin{equation}\label{c-e-p-2}
\langle w'_{(4)}, {\cal Y}_1(w_{(1)}, x_1) {\cal Y}_2(w_{(2)},
x_2)w_{(3)}\rangle_{W_4} \lbar_{x_1= z_1, \;x_2=z_2}
\end{equation}
is absolutely convergent when $|z_1|>|z_2|>0$ and can be analytically
extended to the multivalued analytic function
\begin{equation}\label{c-e-p-3}
\sum_{k=1}^{M}z_2^{r_{k}}(z_1-z_2)^{s_{k}}(\log z_2)^{i_{k}}
(\log(z_1-z_2))^{j_{k}}f_{k}\left(\frac{z_1-z_2}{z_2}\right)
\end{equation}
(here $\log (z_{1}-z_{2})$ and $\log z_{2}$, and in particular, the
powers of the variables, mean the multivalued functions, not the
particular branch we have been using) in the region
$|z_2|>|z_1-z_2|>0$.

\item[Convergence and extension property without logarithms for
products] When $i_{k}=j_{k}=0$ for $k=1, \dots, M$, we call the
property above the {\it convergence and extension property without
logarithms for products}.

\item[Convergence and extension property for iterates] For any
$\beta\in \tilde{A}$, there exists an integer $\tilde{N}_{\beta}$
depending only on ${\cal Y}^1$, ${\cal Y}^2$ and $\beta$, and for any
$w_{(1)}\in W_1$, $w_{(2)}\in (W_2)^{(\beta_{2})}$, $w_{(3)}\in
(W_3)^{(\beta_{3})}$ ($\beta_{2}, \beta_{3} \in \tilde{A}$) and
$w'_{(4)}\in W'_4$, with $w_{(2)}$ and $w_{(3)}$ weight-homogeneous,
such that
\begin{equation}\label{c-e-p-4}
\beta_{2}+\beta_{3}=-\beta, 
\end{equation}
there exist $\tilde M\in{\mathbb N}$, $\tilde{r}_{k}, \tilde{s}_{k}\in
{\mathbb R}$, $\tilde{i}_{k}, \tilde{j}_{k}\in {\mathbb N}$,
$k=1,\dots,\tilde M$; $\tilde{K} \in \Z_{+}$ independent of $w_{(2)}$
and $w_{(3)}$ such that each $\tilde{i}_{k} < \tilde{K}$;
and analytic functions $\tilde{f}_{k}(z)$ on
$|z|<1$, $k=1, \dots, \tilde{M}$, satisfying
\begin{equation}\label{c-e-p-5}
\wt w_{(2)}+\wt w_{(3)}+\tilde{s}_{k}>\tilde{N}_{\beta}, \;\;\;k=1, \dots,
\tilde{M},
\end{equation}
such that
\begin{equation}\label{c-e-p-6}
\langle w'_{(4)}, {\cal Y}^1({\cal Y}^2(w_{(1)}, x_0)w_{(2)},
x_2)w_{(3)}\rangle_{W_4} \lbar_{x_0=z_1-z_2,\;x_2=z_2}
\end{equation}
is absolutely convergent when $|z_2|>|z_1-z_2|>0$ and can be
analytically extended to the multivalued analytic function
\begin{equation}\label{c-e-p-7}
\sum_{k=1}^{\tilde{M}} z_1^{\tilde{r}_{k}}z_2^{\tilde{s}_{k}} (\log
z_1)^{\tilde{i}_{k}} (\log z_2)^{\tilde{j}_{k}}\tilde{f}_{k}
\left(\frac{z_2}{z_1}\right)
\end{equation}
(here $\log z_{1}$ and $\log z_{2}$, and in particular, the powers of
the variables, mean the multivalued functions, not the particular
branches we have been using) in the region $|z_1|>|z_2|>0$.

\item[Convergence and extension property without logarithms 
for iterates] When $i_{k}=j_{k}=0$ for $k=1, \dots, M$, 
we call the property above the {\it convergence and extension 
property without logarithms 
for iterates}.

\end{description}

If the convergence and extension property (with or without logarithms)
for products holds for any objects $W_1$, $W_2$, $W_3$, $W_4$ and
$M_1$ of ${\cal C}$ and any logarithmic intertwining operators ${\cal
Y}_1$ and ${\cal Y}_2$ of the types as above, we say that {\it the
(corresponding) convergence and extension property for products holds
in ${\cal C}$}.  We similarly define the meaning of the phrase {\it
the (corresponding) convergence and extension property for iterates
holds in ${\cal C}$}.

\begin{rema}\label{otherconvergence}
{\rm When we verify these properties using differential equations
below, we will actually be verifying stronger absolute convergence
assertions than as in (\ref{c-e-p-2}) and (\ref{c-e-p-6}) (that is, as
in (\ref{4itm}) and (\ref{4prm})): When all four vectors are
weight-homogeneous, (\ref{c-e-p-2}) and (\ref{c-e-p-6}) will be proved
to be absolutely convergent in the substitution-sense described in the
statement of Proposition \ref{formal=proj}, and this implies the
absolute convergence of (\ref{c-e-p-2}) and (\ref{c-e-p-6}) as above
(again, as in (\ref{4itm}) and (\ref{4prm})), because the absolute
convergence of the triple series (\ref{triple-sum}) implies the
absolute convergence of the corresponding iterated series.}
\end{rema}

\begin{lemma}\label{power-wt}
If the convergence and extension property for products holds,
then, with all four vectors assumed weight-homogeneous, we can always
find $M$, $r_{k}$, $s_{k}$, $i_{k}$, $j_{k}$ and $f_{k}(z)$, $k=1,
\dots, M$, such that
\begin{equation}\label{power-wt-1}
r_{k}+s_{k}=\Delta,\;\;\;k=1, \dots, M,
\end{equation}
where
\begin{equation}
\Delta=-\wt w_{(1)}-\wt w_{(2)} -\wt w_{(3)}+\wt w'_{(4)}.
\end{equation}
Analogously, if the convergence and extension property for iterates
holds, then, with all four vectors assumed weight-homogeneous, we can
always find $\tilde{M}$, $\tilde{r}_{k}$, $\tilde{s}_{k}$,
$\tilde{i}_{k}$, $\tilde{j}_{k}$ and $\tilde{f}_{k}(z)$, $k=1, \dots,
\tilde{M}$, such that
\begin{equation}
\tilde{r}_{k}+\tilde{s}_{k}=\Delta,\;\;\;k=1, \dots, \tilde{M}.
\end{equation}
\end{lemma}
\pf 
Here we prove the assertion for products; the proof for iterates is
entirely analogous.  By (\ref{log:p2}), the 
expression (\ref{c-e-p-2}) but without the evaluation
at $z_1$ and $z_2$ is
\begin{eqnarray*}
\lefteqn{\langle w'_{(4)}, {\cal Y}_1(w_{(1)}, x_1) {\cal Y}_2(w_{(2)},
x_2)w_{(3)}\rangle}\nn
&&=\langle y^{L'(0)}w'_{(4)}, \Y_{1}(y^{-L(0)}w_{(1)}, x_{1}y^{-1})
\Y_{2}(y^{-L(0)}w_{(2)}, x_{2}y^{-1})y^{-L(0)}w_{(3)}\rangle,
\end{eqnarray*}
and we know that this formal series in $x_1$, $\log x_1$, $x_2$, $\log
x_2$, $y$ and $\log y$ is constant with respect to the formal
variables $y$ and $\log y$.  This formal series equals
\begin{eqnarray*}
\lefteqn{\langle y^{\swt w'_{(4)}}e^{\log y (L'(0)-\swt w'_{(4)})}w'_{(4)},
\Y_{1}(y^{-\swt w_{(1)}}e^{-\log y (L(0)-\swt w_{(1)})}w_{(1)}, x_1
y^{-1}) \cdot}\nn
&& \quad
\cdot \Y_{2}(y^{-\swt w_{(2)}}e^{-\log y (L(0)-\swt w_{(2)})}w_{(2)}, x_2 y^{-1})
y^{-\swt w_{(3)}}e^{-\log y (L(0)-\swt w_{(3)})}w_{(3)}\rangle \nn
&& = y^{\Delta}
\langle e^{\log y (L'(0)-\swt w'_{(4)})}w'_{(4)},
\Y_{1}(e^{-\log y (L(0)-\swt w_{(1)})}w_{(1)}, x_1 y^{-1}) \cdot \nn
&& \quad\quad \cdot
\Y_{2}(e^{-\log y (L(0)-\swt w_{(2)})}w_{(2)}, x_2 y^{-1})
e^{-\log y (L(0)-\swt w_{(3)})}w_{(3)}\rangle,
\end{eqnarray*}
and the coefficient of each monomial $x_1^m x_2^n$ ($m,n \in \R$) in
$x_1$ and $x_2$ is a real power of $y$ times a polynomial in $\log
x_1$, $\log x_2$ and $\log y$, and is in fact constant with respect to
$y$ and $\log y$; in particular, $m+n=\Delta$ in each nonzero term.
Thus this expression equals the result of specializing $y$ to $z_2$,
with $\log y$ specialized to $\log z_2$ (using our usual branch of
$\log$), namely,
\begin{eqnarray*}
\lefteqn{e^{\Delta \log z_2}
\langle e^{(\log z_2) (L'(0)-\swt w'_{(4)})}w'_{(4)},
\Y_{1}(e^{-(\log z_2) (L(0)-\swt w_{(1)})}w_{(1)}, e^{-\log z_2} x_1)\cdot}\nn
&& \quad\cdot
\Y_{2}(e^{-(\log z_2) (L(0)-\swt w_{(2)})}w_{(2)}, e^{-\log z_2} x_2)
e^{-(\log z_2) (L(0)-\swt w_{(3)})}w_{(3)}\rangle,
\end{eqnarray*}
using the notation (\ref{log:subs}) for each of $\Y_1$ and $\Y_2$.
Thus, using the notation (\ref{prodabbr}) with $p=0$, we have that
\[
\langle w'_{(4)}, {\cal Y}_1(w_{(1)}, x_1) {\cal Y}_2(w_{(2)},
x_2)w_{(3)}\rangle_{W_4} \lbar_{x_1= z_1, \;x_2=z_2},
\]
an absolutely convergent sum dictated by the monomials in $x_1$ and
$x_2$, equals
\begin{eqnarray*}
\lefteqn{e^{\Delta \log z_2}
\langle e^{(\log z_2) (L'(0)-\swt w'_{(4)})}w'_{(4)},
\Y_{1}(e^{-(\log z_2) (L(0)-\swt w_{(1)})}w_{(1)}, \xi_1)\cdot}\nn
&& \quad\cdot
\Y_{2}(e^{-(\log z_2) (L(0)-\swt w_{(2)})}w_{(2)}, \xi_2)
e^{-(\log z_2) (L(0)-\swt w_{(3)})}w_{(3)}\rangle
\lbar_{\xi_1=e^{\log z_1 - \log z_2}, \;\xi_2=e^0}.
\end{eqnarray*}
Thus (\ref{c-e-p-2}) as an analytic function defined on $|z_{1}|>|z_{2}|>0$ can be 
analytically extended to the analytic function 
\begin{eqnarray}\label{product-1}
\lefteqn{e^{\Delta\log z_{2}}\langle e^{(\log z_{2})(L'(0)-\swt w'_{(4)})}w'_{(4)}, 
\Y_{1}(e^{(\log z_{2})(-L(0)+\swt w_{(1)})}w_{(1)}, 1+(z_{1}-z_{2})/z_{2})\cdot}\nn
&&\quad
\cdot \Y_{2}(e^{(\log z_{2})(-L(0)+\swt w_{(2)})}w_{(2)}, 1)
e^{(\log z_{2})(-L(0)+\swt w_{(3)})}w_{(3)}\rangle
\end{eqnarray}
on the same region, where as usual, we are using the particular branch of $\log$ and the
corresponding branch of each power function.

On the other hand, by
the convergence and extension property for products, (\ref{c-e-p-2})
can be analytically extended to the multivalued function
\begin{eqnarray}\label{product-2}
\lefteqn{\sum_{k=1}^{M}z_2^{r_{k}}(z_1-z_2)^{s_{k}}(\log z_2)^{i_{k}}
(\log(z_1-z_2))^{j_{k}}f_{k}\left(\frac{z_1-z_2}{z_2}\right)}\nn
&&=\sum_{k=1}^{M}z_2^{r_{k}+s_{k}}\left(\frac{z_1-z_2}{z_2}\right)^{s_{k}}
(\log z_2)^{i_{k}}
(\log(z_1-z_2))^{j_{k}}f_{k}\left(\frac{z_1-z_2}{z_2}\right)
\end{eqnarray}
in the region $|z_{2}|>|z_{1}-z_{2}|>0$.  Thus the right-hand side of
(\ref{product-1}) can be analytically extended to the right-hand side
of (\ref{product-2}) in the region $|z_{2}|>|z_{1}-z_{2}|>0$.  Fix
$w\in \C$ such that $|1+w|>1>|w|>0$. Then for any $z_{2}\in
\C^{\times}$, if we let $z_{1}=z_{2}+z_{2}w$, we have
\[
|z_{1}|=|z_{2}||1+w|>|z_{2}|>|z_{2}||w|=|z_{1}-z_{2}|>0,
\]
and thus 
\begin{eqnarray*}
\lefteqn{e^{\Delta\log z_{2}}\langle e^{(\log z_{2})(L'(0)-\swt w'_{(4)})}w'_{(4)}, 
\Y_{1}(e^{(\log z_{2})(-L(0)+\swt w_{(1)})}w_{(1)}, 1+w)\cdot}\nn
&&\quad
\cdot \Y_{2}(e^{(\log z_{2})(-L(0)+\swt w_{(2)})}w_{(2)},
1)e^{(\log z_{2})(-L(0)+\swt w_{(3)})}w_{(3)}\rangle
\end{eqnarray*}
as an analytic function of $z_{2}$ in the region $z_{2}\in
\C^{\times}$ can be analytically extended to the multivalued analytic
function
\[
\sum_{k=1}^{M}z_2^{r_{k}+s_{k}}w^{s_{k}}
(\log z_2)^{i_{k}}
(\log(z_2w))^{j_{k}}f_{k}(w)
\]
in the same region.  In particular, we can change $f_{k}$ if necessary
so that
\begin{eqnarray}\label{product-3}
\lefteqn{e^{\Delta\log z_{2}}\langle e^{(\log z_{2})(L'(0)-\swt w'_{(4)})}
w'_{(4)}, \Y_{1}(e^{(\log z_{2})(-L(0)+\swt w_{(1)})}w_{(1)}, 1+w)\cdot}\nn
&&\quad\quad
\Y_{2}(e^{(\log z_{2})(-L(0)+\swt w_{(2)})}w_{(2)}, 1)
e^{(\log z_{2})(-L(0)+\swt w_{(3)})}w_{(3)}\rangle\nn
&&=\sum_{k=1}^{M}e^{(r_{k}+s_{k})\log z_{2}}w^{s_{k}}
(\log z_2)^{i_{k}}
(\log(z_2w))^{j_{k}}f_{k}(w)
\end{eqnarray}
in the region $z_{2}\in \C^{\times}$.
By Proposition \ref{real-exp-set} (applied to a finite sum), 
we obtain from (\ref{product-3}) that  for $l\ne \Delta$,
\[
\sum_{r_{k}+s_{k}=l}w^{s_{k}}
(\log z_2)^{i_{k}}
(\log(z_2w))^{j_{k}}f_{k}(w)=0
\]
in the region $z_{2}\in \C^{\times}$. Thus we see
that (\ref{c-e-p-2}) can be analytically extended to the multivalued
function (\ref{c-e-p-3}) where $r_{k}+s_{k}=\Delta$.
\epfv

Recall the notion of lower bounded (strongly $\tilde{A}$-graded)
generalized $V$-module in Definition \ref{def:dgw}.  We also need the
following more general notion, for which we recall the notions of
doubly-graded generalized $V$-module and doubly-graded $V$-module in
Definition \ref{doublygraded} (for which the grading restrictions
(\ref{set:dmltc}) and (\ref{set:dmfin}) are not assumed):

\begin{defi}\label{d-g-lower-b}
{\rm If a doubly-graded generalized $V$-module $W=\coprod_{\beta\in
\tilde{A}} \coprod_{n\in {\mathbb R}}W^{(\beta)}_{[n]}$ satisfies the
condition that for $\beta\in \tilde{A}$, $W^{(\beta)}_{[n]}=0$ for $n$
sufficiently negative, we say that $W$ is a {\it lower bounded
doubly-graded generalized $V$-module}.  We define the notion of {\it
lower bounded doubly-graded $V$-module} $W=\coprod_{\beta\in
\tilde{A}} \coprod_{n\in {\mathbb R}}W^{(\beta)}_{(n)}$ analogously.
Such a structure is a lower bounded generalized $V$-module
(respectively, lower bounded $V$-module) if and only if each space
$W^{(\beta)}_{[n]}$ (respectively, $W^{(\beta)}_{(n)}$) is finite
dimensional.}
\end{defi}

We use such lower boundedness to give conditions for insuring that the
convergence and the expansion condition for intertwining maps in
${\cal C}$ hold:

\begin{theo}\label{thm-11.1}
Suppose that the following two conditions are satisfied:
\begin{enumerate}
\item Every finitely-generated lower bounded doubly-graded generalized
$V$-module is an object of ${\cal C}$ (or every finitely-generated
lower bounded doubly-graded $V$-module is an object of ${\cal C}$,
when ${\cal C}$ is in $\mathcal{M}_{sg}$).

\item The convergence and extension property for either products or
iterates holds in ${\cal C}$ (or the convergence and extension
property without logarithms for either products or iterates holds in
${\cal C}$, when ${\cal C}$ is in $\mathcal{M}_{sg}$).
\end{enumerate}
Then the convergence and expansion conditions for intertwining maps in
${\cal C}$ both hold (recall Definitions \ref{conv-conditions} and
\ref{expansion-conditions}).
\end{theo}
\pf By the convergence and extension property for either products or
iterates, the convergence condition for intertwining maps in ${\cal
C}$ holds (recall Definition \ref{conv-conditions}).

We shall now use the convergence and extension property for products
to prove the first of the two equivalent conditions in Theorem
\ref{expansion}; the convergence and extension property for iterates
can be used analogously to prove the second condition in Theorem
\ref{expansion}.  We shall work in the general (logarithmic) case; the
argument for the case ${\cal C}$ in $\mathcal{M}_{sg}$ is analogous
(and shorter).  That is, we shall prove: For any objects $W_1$, $W_2$,
$W_3$, $W_4$ and $M_1$ of ${\cal C}$, any nonzero complex numbers
$z_1$ and $z_2$ satisfying $|z_1|>|z_2|>|z_1-z_2|>0$, any
$P(z_1)$-intertwining map $I_1$ of type ${W_4}\choose {W_1 M_1}$ and
$P(z_2)$-intertwining map $I_2$ of type ${M_1}\choose {W_2W_3}$, and
any $w'_{(4)}\in W'_4$, $(I_1\circ (1_{W_1}\otimes I_2))'(w'_{(4)})\in
(W_1\otimes W_2\otimes W_3)^{*}$ satisfies the
$P^{(2)}(z_1-z_2)$-local grading restriction condition.  Moreover, for
any $w_{(3)}\in W_{3}$ and $n\in \R$, the smallest doubly graded
subspace of $W^{(2)}_{(I_1\circ (1_{W_1}\otimes I_2))'(w'_{(4)}),
w_{(3)}}$ containing the term $\lambda_{n}^{(2)}$ of the (unique)
series $\sum_{n\in \R}\lambda_{n}^{(2)}$ weakly absolutely convergent
to $\mu^{(2)}_{(I_1\circ (1_{W_1}\otimes I_2))'(w'_{(4)}), w_{(3)}}$
as indicated in the $P^{(2)}(z_0)$-grading condition and stable under
the action of $V$ and of $\mathfrak{sl}(2)$ is a generalized
$V$-submodule of some object of ${\cal C}$ included in $(W_{1}\otimes
W_{2})^{*}$.

We may and do take the elements $w_{(1)}\in W_1,w_{(2)}\in W_2,
w_{(3)}\in W_3, w'_{(4)}\in W'_4$ to be weight-homogeneous, since it
is enough to prove the required properties for such elements. We also
recall the notation $\Delta$ from Lemma \ref{power-wt}.  Let
$\Y_{1}=\Y_{I_{1}, 0}$ and $\Y_{2}=\Y_{I_{2}, 0}$.  By Condition 2
(for products) and Lemma \ref{power-wt}, for $\beta\in \tilde{A}$,
there exists $N_{\beta}\in \Z$ and for any $w_{(1)}\in
(W_1)^{(\beta_{1})}_{[n_{1}]}$, $w_{(2)}\in
(W_2)^{(\beta_{2})}_{[n_{2}]}$, $w_{(3)}\in
(W_3)_{[n_{3}]},w'_{(4)}\in (W'_4)_{[n_{4}]}$ such that
(\ref{c-e-p-0}) holds, there exist $M\in{\mathbb N}$, $r_{k}, s_{k}\in
{\mathbb R}$, $i_{k}, j_{k}\in {\mathbb N}$, $k=1,\dots,M$; $K \in
\Z_{+}$ independent of $w_{(1)}$ and $w_{(2)}$ such that each $i_{k} <
K$; and analytic functions $f_{k}(z)$ on $|z|<1$, $k=1, \dots, M$,
such that (\ref{power-wt-1}) holds and (\ref{c-e-p-2}) is absolutely
convergent when $|z_1|>|z_2|>0$ and can be analytically extended to
the multivalued analytic function (\ref{c-e-p-3}) in the region
$|z_2|>|z_1-z_2|>0$.  Then we can always find $f_{k}(z)$ for $k=1,
\dots, M$ such that (\ref{c-e-p-2}) is equal to (\ref{c-e-p-3}) when
we choose the values of $\log z_{2}$ and $\log (z_{1}-z_{2})$ to be
$\log |z_{2}|+i\arg z_{2}$ and $\log |z_{1}-z_{2}|+i\arg
(z_{1}-z_{2})$ where $0\le \arg z_{2}, \arg (z_{1}-z_{2})<2\pi$ and
the values of $z_{2}^{r_{k}}$ and $(z_{1}-z_{2})^{s_{k}}$ to be
$e^{r_{k}\log z_{2}}$ and $e^{s_{k}\log (z_{1}-z_{2})}$ (recalling
that in the convergence and extension properties for products and
iterates, $\log z_{2}$ and $\log (z_{1}-z_{2})$ mean the multivalued
logarithm functions).

Expanding $f_{k}(z)$, $k=1, \dots, M$, we see that in the region
$|z_{1}|>|z_{2}|>|z_{1}-z_{2}|>0$, (\ref{c-e-p-2}) is equal to the
absolutely convergent sum
\begin{eqnarray}\label{thm-11.1-1}
\lefteqn{\sum_{k=1}^{M}
\sum_{m\in \mathbb{N}} C_{m;k}(w'_{(4)},
w_{(1)}, w_{(2)}, w_{(3)})\cdot}\nn &&\quad\quad \cdot
e^{(r_k-m)\log z_{2}}e^{(s_k+m) \log (z_{1}-z_{2})}(\log
z_{2})^{i_k}(\log (z_{1}-z_{2}))^{j_k},
\end{eqnarray}
where, as in this work except in a few identified places, $\log z_{2}$
and $\log (z_{1}-z_{2})$ are the values of the logarithm function at
$z_{2}$ and $z_{1}-z_{2}$ such that $0\le \arg z_{2}, \arg
(z_{1}-z_{2})<2\pi$. 

The expression (\ref{thm-11.1-1}) can be written as 
\begin{eqnarray}\label{thm-11.1-3}
\lefteqn{\sum_{i,j\in\N}\sum_{n\in \mathbb{R}}a_{n; i,j}(w'_{(4)}, 
w_{(1)},
w_{(2)}, w_{(3)})\cdot}\nn
&&\quad\quad
 \cdot e^{(\Delta+n+1)\log (z_{1}-z_{2})}e^{(-n-1)\log z_{2}}
(\log z_{2})^{i}(\log (z_{1}-z_{2}))^{j},
\end{eqnarray} where 
\begin{equation}\label{thm-11.1-5}
a_{n; i,j}(w'_{(4)}, w_{(1)},
w_{(2)}, w_{(3)})=0 
\end{equation}
whenever 
\begin{equation}\label{thm-11.1-4}
n\neq -r_k-1 = -\Delta+s_k-1 \quad\mbox{ or }\quad i\neq i_k \quad\mbox{ or }\quad j\neq j_k
\end{equation}for each $k=1,\cdots,M$.
Also, if $w_{(1)}\in (W_1)^{(\beta_{1})}$, $w_{(2)}\in
(W_2)^{(\beta_{2})}$, $w_{(3)}\in (W_3)^{(\beta_{3})}$ and
$w'_{(4)}\in (W'_4)^{(\beta_{4})}$ and
\begin{equation}\label{beta-sum}
\beta_{1}+\beta_{2}+\beta_{3}+\beta_{4}\ne 0,
\end{equation}
then (\ref{thm-11.1-5}) holds.
{}From (\ref{thm-11.1-4}), (\ref{beta-sum}), (\ref{c-e-p-1}), (\ref{power-wt-1})
we see that for $n\in \mathbb{R}$, if 
\begin{equation}\label{thm-11.1-6}
n+1+ \wt w'_{(4)}-
\wt w_{(3)}\le N_{\beta_{3}+\beta_{4}},
\end{equation}
then (\ref{thm-11.1-5}) holds.

Since $|z_{1}|>|z_{2}|>|z_{1}-z_{2}|>0$, we know that
\begin{eqnarray}\label{thm-11.1-7}
\lefteqn{\sum_{i,j\in\N}\sum_{n\in \mathbb{R}}a_{n; i,j}(w'_{(4)}, 
w_{(1)},
w_{(2)}, w_{(3)})\cdot}\nn
&&\quad
 \cdot e^{(\Delta+n+1)\log (z_{1}-z_{2})}e^{(-n-1)\log z_{2}}
(\log z_{2})^{i}(\log (z_{1}-z_{2}))^{j}
\end{eqnarray}
converges absolutely to (\ref{c-e-p-2}).  For $i,j\in\N$, $w'_{(4)}\in W'_{4}$,
$w_{(3)}\in W_{3}$, let $\beta_{n; i,j}(w'_{(4)}, w_{(3)})\in
(W_{1}\otimes W_{2})^{*}$ be defined by
\begin{eqnarray}\label{thm-11.1-8}
\lefteqn{(\beta_{n; i,j}(w'_{(4)}, w_{(3)}))(w_{(1)}\otimes w_{(2)})}\nn
&&=a_{n; i,j}(e^{L(0)_s\log(z_1-z_2)}w'_{(4)},
e^{-L(0)_s\log(z_1-z_2)}w_{(1)},e^{-L(0)_s\log(z_1-z_2)}w_{(2)},
e^{-L(0)_s\log(z_1-z_2)}w_{(3)})
\nn\end{eqnarray}
for all $w_{(1)}\in W_{1}$ and $w_{(2)}\in W_{2}$.
By definition, the series 
\[
\sum_{i,j\in\N}\sum_{n\in \mathbb{R}}(\beta_{n; i,j}
(w'_{(4)}, w_{(3)}))(w_{(1)}\otimes w_{(2)})
e^{(n+1)\log (z_{1}-z_{2})}e^{(-n-1)\log z_{2}}
(\log z_{2})^{i}(\log (z_{1}-z_{2}))^{j}
\]
is absolutely 
convergent to 
$\mu^{(2)}_{I_{1}\circ (1_{W_{2}}\otimes I_{2}))'(w'_{(4)}), w_{(3)}}
(w_{(1)}\otimes w_{(2)})$ for $w_{(1)}\in W_{1}$ and $w_{(2)}\in W_{2}$.
Also since (\ref{thm-11.1-5}) holds when (\ref{beta-sum}) holds, 
we have
\begin{equation}\label{beta-Atilde}
\beta_{n; i,j}(w'_{(4)}, w_{(3)})\in ((W_{1}\otimes W_{2})^{*})^{(\beta_{3}+\beta_{4})}
\end{equation}
for $w_{(3)}\in (W_3)^{(\beta_{3})}$
and $w'_{(4)}\in (W'_4)^{(\beta_{4})}$.

To show that $(I_{1}\circ (1_{W_{2}}\otimes I_{2}))'(w'_{(4)})$
satisfies the $P^{(2)}(z_{1}-z_{2})$-local grading restriction
condition, we need to calculate
\[
(v_{1})_{m_{1}}\cdots (v_{r})_{m_{r}}\beta_{n; i,j}(w'_{(4)}, 
w_{(3)})
\]
and its weight for any $r\in \mathbb{N}$, $v_{1}, \dots, v_{r}\in V,$
$m_{1}, \dots, m_{r}\in \mathbb{Z}$, $n\in \mathbb{R}$, where
$(v_{1})_{m_{1}}, \dots, (v_{r})_{m_{r}}$, $m_{1}, \cdots, m_{r}\in
\mathbb{Z}$, are the components of $Y'_{P(z_{1}-z_{2})}(v_{1}, x)$,
$\dots\ $, $Y'_{P(z_{1}-z_{2})}(v_{r}, x)$, respectively, on
$(W_{1}\otimes W_{2})^{*}$. As in \cite{tensor4}, for convenience, we
instead calculate
\[
(v^{o}_{1})_{m_{1}}\cdots (v^{o}_{r})_{m_{r}}
\beta_{n; i,j}(w'_{(4)}, w_{(3)}),
\]
where $(v^{o}_{1})_{m_{1}}, \dots, (v^{o}_{r})_{m_{r}}$, $m_{1},
\cdots, m_{r}\in \mathbb{Z}$, are the components of the opposite
vertex operators $Y^{\prime\; o}_{P(z_{1}-z_{2})}(v_{1}, x)$, $\dots\
$, $Y^{\prime\; o}_{P(z_{1}-z_{2})}(v_{r}, x)$, respectively.

By the definition (\ref{Y'def}) of
$Y'_{P(z_{1}-z_{2})}(v, x)$ and 
\[
Y^{\prime\; o}_{P(z_{1}-z_{2})}(v, x)=Y'_{P(z_{1}-z_{2})}
(e^{xL(1)}(-x^{-2})^{L(0)}v, x^{-1}),
\]
we have
\begin{eqnarray}\label{16.19}
\lefteqn{(Y^{\prime \; o}_{P(z_{1}-z_{2})}(v, x)
\beta_{n; i,j}(w'_{(4)}, w_{(3)}))(w_{(1)}\otimes w_{(2)})}\nno\\
&&=(\beta_{n; i,j}(w'_{(4)}, w_{(3)}))(w_{(1)}\otimes Y_{2}(v, x)
w_{(2)})
\nno\\
&&\quad +\res_{x_{0}}(z_{1}-z_{2})^{-1}\delta\left(\frac{x-x_{0}}{z_{1}-z_{2}}
\right)
(\beta_{n; i,j}(w'_{(4)}, w_{(3)}))(Y_{1}(v, x_{0})
w_{(1)}\otimes w_{(2)})
\nno\\
&&=a_{n; i,j}(e^{L(0)_s\log(z_1-z_2)}w'_{(4)},
e^{-L(0)_s\log(z_1-z_2)}w_{(1)}, e^{-L(0)_s\log(z_1-z_2)}Y_{2}(v, x)
w_{(2)}, e^{-L(0)_s\log(z_1-z_2)}w_{(3)})\nno\\
&&\quad +\res_{x_{0}}(z_{1}-z_{2})^{-1}\delta\left(\frac{x-x_{0}}{z_{1}-z_{2}}
\right)\cdot\nn
&&\cdot a_{n; i,j}(e^{L(0)_s\log(z_1-z_2)}w'_{(4)}, e^{-L(0)_s\log(z_1-z_2)}Y_{1}(v, x_{0})
w_{(1)}, e^{-L(0)_s\log(z_1-z_2)}w_{(2)}, e^{-L(0)_s\log(z_1-z_2)}w_{(3)}).\nn
\end{eqnarray}
On the other hand, since (\ref{thm-11.1-7}) is absolutely 
convergent to 
(\ref{c-e-p-2})
when $|z_{1}|>|z_{2}|>|z_{1}-z_{2}|>0$, we have
\begin{eqnarray}\label{16.20}
\lefteqn{\sum_{i,j\in\N}\sum_{n\in\R} e^{(n+1)\log
(z_{1}-z_{2})}e^{(-n-1)\log z_{2}}(\log z_{2})^{i}(\log
(z_{1}-z_{2}))^{j}\cdot
}\nno\\
&& \cdot a_{n; i,j}(e^{L(0)_s\log(z_1-z_2)}w'_{(4)}, e^{-L(0)_s\log(z_1-z_2)}w_{(1)},
e^{-L(0)_s\log(z_1-z_2)}Y_{2}(v, x)w_{(2)}, e^{-L(0)_s\log(z_1-z_2)}w_{(3)}) \nno\\
&& +\res_{x_{0}}(z_{1}-z_{2})^{-1}\delta\left(\frac{x-x_{0}}{z_{1}-z_{2}}
\right)\sum_{i,j\in\N}\sum_{n\in \mathbb{R}} 
e^{(n+1)\log (z_{1}-z_{2})}e^{(-n-1)\log z_{2}}
(\log z_{2})^{i}(\log (z_{1}-z_{2}))^{j}\cdot\nno\\
&&\cdot a_{n; i,j}(e^{L(0)_s\log(z_1-z_2)}w'_{(4)},e^{-L(0)_s\log(z_1-z_2)} Y_{1}(v, x_{0})
w_{(1)},e^{-L(0)_s\log(z_1-z_2)} w_{(2)},e^{-L(0)_s\log(z_1-z_2)} w_{(3)})\nno\\
&&=\langle w'_{(4)}, {\cal Y}_{1}(w_{(1)}, x_{1})
{\cal Y}_{2}(Y_{2}(v, x)w_{(2)}, x_{2})w_{(3)}\rangle_{W_{4}}
\lbar_{x_{1}
=z_{1}, x_{2}=z_{2}}\nno\\
&&\quad +\res_{x_{0}}(x_{1}-x_{2})^{-1}\delta\left(\frac{x-x_{0}}{x_{1}-x_{2}}
\right)\cdot\nno\\
&&\quad\quad\quad\quad\cdot
\langle w'_{(4)}, {\cal Y}_{1}(Y_{1}(v, x_{0})w_{(1)},
x_{1})
{\cal Y}_{2}(w_{(2)}, x_{2})w_{(3)}\rangle_{W_{4}}
\lbar_{x_{1}
=z_{1}, x_{2}=z_{2}}\nno\\
&&=\langle w'_{(4)}, {\cal Y}_{1}(w_{(1)}, x_{1})
{\cal Y}_{2}(Y_{2}(v, x)w_{(2)}, x_{2})w_{(3)}\rangle_{W_{4}}
\lbar_{x_{1}
=z_{1}, x_{2}=z_{2}}\nno\\
&&\quad +\res_{x_{0}}x_{1}^{-1}\delta\left(\frac{(x+x_{2})-x_{0}}{x_{1}}
\right)\cdot\nno\\
&&\quad\quad\quad\quad\cdot
\langle w'_{(4)}, {\cal Y}_{1}(Y_{1}(v, x_{0})w_{(1)},
x_{1})
{\cal Y}_{2}(w_{(2)}, x_{2})w_{(3)}\rangle_{W_{4}}
\lbar_{x_{1}
=z_{1}, x_{2}=z_{2}}.
\end{eqnarray}

Using the Jacobi identity for the logarithmic intertwining operators
${\cal Y}_{1}$ and ${\cal Y}_{2}$ and the properties of the formal
$\delta$-function, the right-hand side of (\ref{16.20}) is equal to
\begin{eqnarray}\label{16.21}
\lefteqn{\res_{y}x^{-1}\delta\left(\frac{y-x_{2}}{x}\right)
\langle w'_{(4)}, {\cal Y}_{1}(w_{(1)}, x_{1})Y_{5}(v, y)
{\cal Y}_{2}(w_{(2)}, x_{2})w_{(3)}\rangle_{W_{4}}
\lbar_{x_{1}
=z_{1}, x_{2}=z_{2}}}\nno\\
&&- \res_{y}x^{-1}\delta\left(\frac{x_{2}-y}{-x}\right)
\langle w'_{(4)}, {\cal Y}_{1}(w_{(1)}, x_{1})
{\cal Y}_{2}(w_{(2)}, x_{2})Y_{3}(v, y)w_{(3)}\rangle_{W_{4}}
\lbar_{x_{1}
=z_{1}, x_{2}=z_{2}}\nno\\
&&+\langle w'_{(4)}, Y_{4}(v, x+x_{2}){\cal Y}_{1}(w_{(1)},
x_{1})
{\cal Y}_{2}(w_{(2)}, x_{2})w_{(3)}\rangle_{W_{4}}
\lbar_{x_{1}
=z_{1}, x_{2}=z_{2}}
\nno\\
&&-\langle w'_{(4)}, {\cal Y}_{1}(w_{(1)},
x_{1})Y_{5}(v, x+x_{2})
{\cal Y}_{2}(w_{(2)}, x_{2})w_{(3)}\rangle_{W_{4}}
\lbar_{x_{1}
=z_{1}, x_{2}=z_{2}}
\nno\\
&&=\langle w'_{(4)}, Y_{4}(v, x+x_{2}){\cal Y}_{1}(w_{(1)},
x_{1})
{\cal Y}_{2}(w_{(2)}, x_{2})w_{(3)}\rangle_{W_{4}}
\lbar_{x_{1}
=z_{1}, x_{2}=z_{2}}\nno\\
&&\quad - \res_{y}x^{-1}\delta\left(\frac{x_{2}-y}{-x}\right)
\langle w'_{(4)}, {\cal Y}_{1}(w_{(1)}, x_{1})
{\cal Y}_{2}(w_{(2)}, x_{2})Y_{3}(v, y)w_{(3)}\rangle_{W_{4}}
\lbar_{x_{1}
=z_{1}, x_{2}=z_{2}}\nno\\
&&=\res_{y}x^{-1}\delta\left(\frac{y-x_{2}}{x}\right)
\langle w'_{(4)}, Y_{4}(v, y){\cal Y}_{1}(w_{(1)},
x_{1})
{\cal Y}_{2}(w_{(2)}, x_{2})w_{(3)}\rangle_{W_{4}}
\lbar_{x_{1}
=z_{1}, x_{2}=z_{2}}\nno\\
&&\quad - \res_{y}x^{-1}\delta\left(\frac{x_{2}-y}{-x}\right)
\langle w'_{(4)}, {\cal Y}_{1}(w_{(1)}, x_{1})
{\cal Y}_{2}(w_{(2)}, x_{2})Y_{3}(v, y)w_{(3)}\rangle_{W_{4}}
\lbar_{x_{1}
=z_{1}, x_{2}=z_{2}}.\nno\\
\end{eqnarray}
{}From (\ref{16.20}) and (\ref{16.21}), we obtain
\begin{eqnarray}\label{16.22}
\lefteqn{\sum_{i,j\in\N}\sum_{n\in\R} e^{(n+1)\log
(z_{1}-z_{2})}e^{(-n-1)\log z_{2}}(\log z_{2})^{i}(\log
(z_{1}-z_{2}))^{j}\cdot
}\nno\\
&& \cdot a_{n; i,j}(e^{L(0)_s\log(z_1-z_2)}w'_{(4)}, e^{-L(0)_s\log(z_1-z_2)}w_{(1)},
e^{-L(0)_s\log(z_1-z_2)}Y_{2}(v, x)w_{(2)}, e^{-L(0)_s\log(z_1-z_2)}w_{(3)}) \nno\\
&& +\res_{x_{0}}(z_{1}-z_{2})^{-1}\delta\left(\frac{x-x_{0}}{z_{1}-z_{2}}
\right)\sum_{i,j\in\N}\sum_{n\in \mathbb{R}} 
e^{(n+1)\log (z_{1}-z_{2})}e^{(-n-1)\log z_{2}}
(\log z_{2})^{i}(\log (z_{1}-z_{2}))^{j}\cdot\nno\\
&&\cdot a_{n; i,j}(e^{L(0)_s\log(z_1-z_2)}w'_{(4)},e^{-L(0)_s\log(z_1-z_2)} Y_{1}(v, x_{0})
w_{(1)},e^{-L(0)_s\log(z_1-z_2)} w_{(2)},e^{-L(0)_s\log(z_1-z_2)} w_{(3)})\nno\\
&&=\res_{y}x^{-1}\delta\left(\frac{y-x_{2}}{x}\right)
\langle w'_{(4)}, Y(v, y){\cal Y}_{1}(w_{(1)}, x_{1})
{\cal Y}_{2}(w_{(2)}, x_{2})w_{(3)}\rangle_{W_{4}}
\lbar_{x_{1}
=z_{1}, x_{2}=z_{2}}\nno\\
&&\quad - \res_{y}x^{-1}\delta\left(\frac{x_{2}-y}{x}\right)
\langle w'_{(4)}, {\cal Y}_{1}(w_{(1)}, x_{1})
{\cal Y}_{2}(w_{(2)}, x_{2})Y(v, y)w_{(3)}\rangle_{W_{4}}
\lbar_{x_{1}
=z_{1}, x_{2}=z_{2}}\nno\\
&&=\sum_{m\in \mathbb{Z}}\sum_{l\in \N}
(-1)^{l}{{m}\choose {l}}x^{-m-1}x_{2}^{l}\langle w'_{(4)}, v_{m-l}
{\cal Y}_{1}(w_{(1)}, x_{1})
{\cal Y}_{2}(w_{(2)}, x_{2})w_{(3)}\rangle_{W_{4}}
\lbar_{x_{1}
=z_{1}, x_{2}=z_{2}}\nno\\
&&\quad -\sum_{m\in \mathbb{Z}}\sum_{ l\in \N}
(-1)^{l+m}{{m}\choose {l}}x^{-m-1}x_{2}^{m-l}\langle w'_{(4)},
{\cal Y}_{1}(w_{(1)}, x_{1})
{\cal Y}_{2}(w_{(2)}, x_{2})v_{l}
w_{(3)}\rangle_{W_{4}}\lbar_{x_{1}
=z_{1}, x_{2}=z_{2}}\nno\\
&&=\sum_{m\in \mathbb{Z}}\sum_{l\in \N}
(-1)^{l}{{m}\choose {l}}x^{-m-1}\sum_{i,j\in\N}
\sum_{n\in \mathbb{R}}e^{(n+1)\log (z_{1}-z_{2})}e^{(l-n-1)\log z_{2}}
(\log z_{2})^{i}(\log (z_{1}-z_{2}))^{j}\cdot\nno\\
&&\cdot a_{n; i,j}(e^{L(0)_s\log(z_1-z_2)}v^{o}_{m-l}w'_{(4)},
e^{-L(0)_s\log(z_1-z_2)}w_{(1)}, e^{-L(0)_s\log(z_1-z_2)}w_{(2)},
e^{-L(0)_s\log(z_1-z_2)}w_{(3)})\nno\\
&&\quad -\sum_{m\in \mathbb{Z}}\sum_{ l\in \N}
(-1)^{l+m}{{m}\choose {l}}x^{-m-1}
\sum_{i,j\in\N}\sum_{n\in \mathbb{R}}e^{(n+1)\log (z_{1}-z_{2})}e^{(m-l-n-1)\log z_{2}}
(\log z_{2})^{i}(\log (z_{1}-z_{2}))^{j}\cdot\nno\\
&&\cdot a_{n; i,j}(e^{L(0)_s\log(z_1-z_2)}w'_{(4)},
e^{-L(0)_s\log(z_1-z_2)}w_{(1)},e^{-L(0)_s\log(z_1-z_2)}w_{(2)},e^{-L(0)_s
\log(z_1-z_2)}v_{l}w_{(3)}).\nn
\end{eqnarray}
The intermediate steps in the equality (\ref{16.22}) hold only when
$|z_{1}|>|z_{2}|>|z_{1}-z_{2}|>0$. But since both sides of
(\ref{16.22}) can be analytically extended to the region
$|z_{2}|>|z_{1}-z_{2}|>0$, (\ref{16.22}) must hold when
$|z_{2}|>|z_{1}-z_{2}|>0$. By Proposition \ref{real-exp-set}, the
coefficients of both sides of
(\ref{16.22}) in powers of $e^{\log z_{2}}$, $\log z_{2}$ and $\log
(z_{1}-z_{2})$ are equal, that is,
\begin{eqnarray}\label{16.23}
\lefteqn{a_{n; i,j}(e^{L(0)_s\log(z_1-z_2)}w'_{(4)}, 
e^{-L(0)_s\log(z_1-z_2)}w_{(1)}, e^{-L(0)_s\log(z_1-z_2)}Y_{2}(v,
x)w_{(2)}, e^{-L(0)_s\log(z_1-z_2)}w_{(3)})}\nno\\
&&+\res_{y}(z_{1}-z_{2})^{-1}\delta\left(\frac{x-y}{z_{1}-z_{2}}
\right)\cdot\nno\\
&&\cdot a_{n; i,j}(e^{L(0)_s\log(z_1-z_2)}w'_{(4)}, e^{-L(0)_s\log(z_1-z_2)}Y_{1}(v, y)
w_{(1)}, e^{-L(0)_s\log(z_1-z_2)}w_{(2)},e^{-L(0)_s\log(z_1-z_2)}w_{(3)})\nno\\
&&=\sum_{m\in \mathbb{Z}}\sum_{l\in \N}
(-1)^{l}{{m}\choose {l}}x^{-m-1}(z_{1}-z_{2})^{l}\cdot\nno\\
&&\cdot a_{n+l; i,j}(e^{L(0)_s\log(z_1-z_2)}v^{o}_{m-l}w'_{(4)},
e^{-L(0)_s\log(z_1-z_2)}w_{(1)}, e^{-L(0)_s\log(z_1-z_2)}w_{(2)},
e^{-L(0)_s\log(z_1-z_2)}w_{(3)})\nno\\
&&\quad -\sum_{m\in \mathbb{Z}}\sum_{ l\in \N}
(-1)^{l+m}{{m}\choose {l}}x^{-m-1}(z_{1}-z_{2})^{m-l}\cdot \nno\\
&&\cdot a_{n+m-l; i,j}(e^{L(0)_s\log(z_1-z_2)}w'_{(4)},
e^{-L(0)_s\log(z_1-z_2)}w_{(1)}, e^{-L(0)_s\log(z_1-z_2)}w_{(2)},
e^{-L(0)_s\log(z_1-z_2)}v_{l}w_{(3)}).\nno\\
&&
\end{eqnarray}

By (\ref{16.19}) and (\ref{16.23}), we obtain
\begin{eqnarray}
\lefteqn{(v^{o}_{m}\beta_{n; i,j}(w'_{(4)}, w_{(3)}))(w_{(1)}
\otimes w_{(2)})}\nno\\
&&=\sum_{l\in \N}
(-1)^{l}{{m}\choose {l}}(\beta_{n+l; i,j}(v^{o}_{m-l}w'_{(4)}, w_{(3)}))
(w_{(1)}, w_{(2)})(z_{1}-z_{2})^{l}\nno\\
&&\quad -\sum_{ l\in \N}
(-1)^{l+m}{{m}\choose {l}}(\beta_{n+m-l; i,j}(w'_{(4)}, v_{l}w_{(3)}))(w_{(1)},
w_{(2)})(z_{1}-z_{2})^{m-l}.
\end{eqnarray}
By induction, we obtain
\begin{eqnarray}\label{thm-11.1-9}
\lefteqn{((v^{o}_{1})_{m_{1}}\cdots (v^{o}_{r})_{m_{r}}\beta_{n; i,j}(w'_{(4)}, 
w_{(3)})
(w_{(1)}\otimes w_{(2)})=}\nno\\
&&=\sum_{s=0}^{t}\sum_{\begin{array}{c}
\mbox{\scriptsize $j_{1}>\cdots >j_{s}$}\\
\mbox{\scriptsize $j_{s+1}>\cdots >j_{r}$}
\\ 
\mbox{\scriptsize $\{j_{1}, \dots, j_{r}\}=\{1, \dots, r\}$}\end{array}}
\sum_{l_{1},  \dots, l_{r}\ge 0}(-1)^{l_{1}+\cdots +l_{r}+(m_{j_{s+1}}+1)+
\cdots
+ (m_{j_{r}}+1)}
{{m_{j_{1}}}\choose {l_{1}}}\cdots {{m_{j_{r}}}\choose {l_{r}}}\cdot\nno\\
&&\hspace{2em}\cdot
(\beta_{n+m_{j_{s+1}}\cdots +m_{j_{r}}+l_{1}+\cdots +l_{s}-l_{s+1}-
\cdots -l_{r}; i,j}\nno\\
&&\hspace{4em}
((v^{o}_{j_{1}})_{m_{j_{1}}-l_{1}}\cdots (v^{o}_{j_{s}})_{m_{j_{s}}-l_{s}}
w'_{(4)}, v_{l_{s+1}}\cdots v_{l_{r}}w_{(3)}))(w_{(1)}\otimes w_{(2)})
\cdot \nno\\
&&\hspace{2em}\cdot (z_{1}-z_{2})^{l_{1}+\cdots +l_{s}+(m_{s+1}-l_{s+1})+\cdots 
+(m_{r}-l_{r})}
\end{eqnarray}
for $m_{1}, \dots, m_{r}\in \mathbb{Z}$,
and any $v_{1}, \dots, 
v_{r}\in V$. In particular, when $V$ is a conformal vertex algebra,
\begin{eqnarray}\label{thm-11.1-10}
\lefteqn{(L'_{P(z_{1}-z_{2})}(0)\beta_{n; i,j}(w'_{(4)}, 
w_{(3)}))(w_{(1)}
\otimes w_{(2)})}\nno\\
&&=(\beta_{n; i,j}(L'(0)w'_{(4)}, w_{(3)}))(w_{(1)}
\otimes w_{(2)})\nno\\
&&\quad -(\beta_{n+1; i,j}(L'(1)w'_{(4)}, w_{(3)}))(w_{(1)}
\otimes w_{(2)})\cdot(z_1-z_2)\nno\\
&&\quad -(\beta_{n; i,j}(w'_{(4)}, L(0)w_{(3)}))(w_{(1)}
\otimes w_{(2)})\nno\\
&&\quad +(\beta_{n+1; i,j}(w'_{(4)}, L(-1)w_{(3)}))(w_{(1)}
\otimes w_{(2)})\cdot(z_1-z_2)\nno\\
&&=(\wt w'_{(4)}-\wt w_{(3)})(\beta_{n; i,j}(w'_{(4)}, w_{(3)}))(w_{(1)}
\otimes w_{(2)})\nno\\
&&\quad +(\beta_{n; i,j}((L'(0)-\wt w'_{(4)})w'_{(4)}, w_{(3)}))(w_{(1)}
\otimes w_{(2)})\nno\\
&&\quad -(\beta_{n; i,j}(w'_{(4)}, (L(0)-\wt w_{(3)})w_{(3)}))(w_{(1)}
\otimes w_{(2)})\nno\\
&&\quad -(\beta_{n+1; i,j}(L'(1)w'_{(4)}, w_{(3)}))(w_{(1)}
\otimes w_{(2)})\cdot(z_1-z_2)\nno\\
&&\quad +(\beta_{n+1; i,j}(w'_{(4)}, L(-1)w_{(3)}))(w_{(1)}
\otimes w_{(2)})\cdot(z_1-z_2);
\end{eqnarray}
by a similar argument, (\ref{thm-11.1-10}) holds for a M\"{o}bius
vertex algebra as well.  {}From (\ref{thm-11.1-9}) and
(\ref{beta-Atilde}), we see that when $v_{1}\in
V_{(m_{1})}^{(\alpha_{1})}, \dots, v_{r}\in
V_{(m_{r})}^{(\alpha_{r})}$, $w_{(3)}\in
(W_{3})_{[n_{3}]}^{(\beta_{3})}$ and $w'_{(4)}\in
(W_{4}')_{[n_{4}]}^{(\beta_{4})}$,
\begin{equation}\label{thm-11.1-10.1}
(v^{o}_{1})_{m_{1}}\cdots (v^{o}_{r})_{m_{r}}\beta_{n; i,j}(w'_{(4)}, 
w_{(3)})\in ((W_{1}\otimes W_{2})^{*})^{(\alpha_{1}+\cdots +\alpha_{r}+\beta_{3}+\beta_{4})}.
\end{equation}

Let $z_{0}=z_{1}-z_{2}$. When $|z_{1}|>|z_{2}|>|z_{0}|>0$, 
\begin{eqnarray}\label{16.27}
\lefteqn{-\sum_{i,j\in\N}\sum_{n\in \mathbb{R}}a_{n; i,j}(L'(1)w'_{(4)},
w_{(1)}, w_{(2)}, w_{(3)})\cdot}\nno\\
&&\quad\quad\quad\quad\cdot 
e^{((\Delta-1)+n+1)\log (z_{1}-z_{2})}e^{(-n-1)\log z_{2}}
(\log z_{2})^{i}(\log (z_{1}-z_{2}))^{j}
\nno\\
&&\quad +\sum_{i,j\in\N}\sum_{n\in \mathbb{R}}a_{n; i,j}(w'_{(4)},
w_{(1)}, w_{(2)}, L(-1)w_{(3)}) \cdot \nno\\
&&\quad\quad\quad\quad\cdot e^{((\Delta-1)+n+1)\log (z_{1}-z_{2})}
e^{(-n-1)\log z_{2}}
(\log z_{2})^{i}(\log (z_{1}-z_{2}))^{j}\nno\\
&&= -
\langle L'(1)w'_{(4)},
{\cal Y}_{1}(
w_{(1)}, x_{1})
{\cal Y}_{2}(w_{(2)}, x_{2})w_{(3)}\rangle_{W_{4}}
\lbar_{x_{1}
=z_{2}+z_{0}, \; x_{2}=z_{2}}\nno\\
&&\quad +
\langle w'_{(4)}, {\cal Y}_{1}(
w_{(1)}, x_{1})
{\cal Y}_{2}(w_{(2)}, x_{2})L(-1)w_{(3)}\rangle_{W_{4}}
\lbar_{x_{1}
=z_{2}+z_{0}, \; x_{2}=z_{2}}.
\end{eqnarray}
By the commutator formula for $L(-1)$ and intertwining operators and
the $L(-1)$-derivative property for intertwining operators, the
right-hand side of (\ref{16.27}) is equal to
\begin{eqnarray}\label{16.28}
\lefteqn{-
\left\langle w'_{(4)},
\left(\frac{d}{dx_{1}}({\cal Y}_{1}(
w_{(1)}, x_{1}))\right)
{\cal Y}_{2}(w_{(2)}, x_{2})w_{(3)})\right\rangle_{W_{4}}
\lbar_{x_{1}
=z_{2}+z_{0}, \; x_{2}=z_{2}}}\nno\\
&&\quad -
\left\langle w'_{(4)}, {\cal Y}_{1}(
w_{(1)}, x_{1})
\left(\frac{d}{dx_{2}}({\cal Y}_{2}(w_{(2)}, x_{2})\right)w_{(3)})\right\rangle_{W_{4}}
\lbar_{x_{1}
=z_{2}+z_{0}, \; x_{2}=z_{2}}\nno\\
&&=-\frac{\p}{\p z_{2}}\biggl(
\langle w'_{(4)},
({\cal Y}_{1}(
w_{(1)}, x_{1})
{\cal Y}_{2}(w_{(2)}, x_{2})w_{(3)})\rangle_{W_{4}}
\lbar_{x_{1}
=z_{2}+z_{0}, \; x_{2}=z_{2}}\biggr)\nno\\
&&=-\frac{\p}{\p z_{2}}\Biggl(\sum_{i,j\in\N}
\sum_{n\in \mathbb{R}}a_{n; i,j}(w'_{(4)},
w_{(1)}, w_{(2)}, w_{(3)}) \cdot \nno\\
&&\quad\quad\quad\quad\quad\quad\quad\quad\cdot
e^{(\Delta+n+1)\log z_{0}}e^{(-n-1)\log z_{2}}
(\log z_{2})^{i}(\log z_{0})^{j}
\Biggr)\nno\\
&&=-\sum_{i,j\in\N}\sum_{n\in \mathbb{R}}(-n-1)a_{n; i,j}(w'_{(4)},
w_{(1)}, w_{(2)}, w_{(3)}) \cdot \nno\\
&&\quad\quad\quad\quad\quad\quad\quad\quad
\cdot e^{(\Delta+n+1)\log z_{0}}e^{(-n-2)\log z_{2}}
(\log z_{2})^{i}(\log z_{0})^{j}\nno\\
&&\quad -\sum_{i,j\in\N}\sum_{n\in \mathbb{R}}a_{n; i,j}(w'_{(4)},
w_{(1)}, w_{(2)}, w_{(3)}) \cdot \nno\\
&&\quad\quad\quad\quad\quad\quad\quad\quad
\cdot e^{(\Delta+n+1)\log z_{0}}e^{(-n-1)\log z_{2}}
\left(\frac{\partial}{\partial z_{2}}
(\log z_{2})^{i}\right)(\log z_{0})^{j},
\end{eqnarray}
where $\frac{\p}{\p z_{2}}$ is partial differentiation with respect to
$z_{2}$ acting on functions of $z_{0}$ and $z_{2}$ rather than on
functions of $z_{1}$ and $z_{2}$.

{}From (\ref{16.27}) and (\ref{16.28}), we obtain
\begin{eqnarray}\label{thm-11.1-11}
\lefteqn{-\sum_{i,j\in\N}\sum_{n\in \mathbb{R}}\beta_{n; i,j}(L'(1)w'_{(4)},
w_{(3)})(w_{(1)}\otimes w_{(2)})\cdot}\nn
&&\quad\quad\quad\quad\quad\quad
 \cdot e^{(n+1)\log (z_{1}-z_{2})}e^{(-n-1)\log z_{2}}
(\log z_{2})^{i}(\log (z_{1}-z_{2}))^{j}
\nno\\
&&\quad +\sum_{i,j\in\N}\sum_{n\in \mathbb{R}}\beta_{n; i,j}(w'_{(4)},
L(-1)w_{(3)})(w_{(1)}\otimes w_{(2)}) \cdot\nn
&&\quad\quad\quad\quad\quad\quad
\cdot e^{(n+1)\log (z_{1}-z_{2})}
e^{(-n-1)\log z_{2}}
(\log z_{2})^{i}(\log (z_{1}-z_{2}))^{j}\nno\\
&&=-\sum_{i,j\in\N}\sum_{n\in \mathbb{R}}(-n-1)\beta_{n; i,j}(w'_{(4)},
w_{(3)})(w_{(1)}\otimes w_{(2)})  \cdot\nn
&&\quad\quad\quad\quad\quad\quad
e^{(n+1)\log (z_{1}-z_{2})}e^{(-n-2)\log z_{2}}
(\log z_{2})^{i}(\log (z_{1}-z_{2}))^{j}
\nno\\
&&\quad-\sum_{i,j\in\N}\sum_{n\in \mathbb{R}}\beta_{n; i,j}(w'_{(4)},
w_{(3)})(w_{(1)}\otimes w_{(2)})  \cdot\nn
&&\quad\quad\quad\quad\quad\quad
e^{(n+1)\log (z_{1}-z_{2})}e^{(-n-1)\log z_{2}}
\left(\frac{\partial}{\partial z_{2}}
(\log z_{2})^{i}\right)(\log (z_{1}-z_{2}))^{j}\nn
\end{eqnarray}
when $|z_{1}|>|z_{2}|>|z_{1}-z_{2}|>0$. 
Since both sides of (\ref{thm-11.1-11}) can be analytically 
extended to the region $|z_{2}|>|z_{1}-z_{2}|>0$, it 
also holds when $|z_{2}|>|z_{1}-z_{2}|>0$. Thus the 
expansion coefficients of the two sides of (\ref{thm-11.1-11})
as series in $z_{2}$ and $z_{1}-z_{2}$ are equal. So 
for $i,j\in\N$,
\begin{eqnarray}\label{thm-11.1-12}
\lefteqn{(z_1-z_2)\Big(-\beta_{n+1; i,j}(L'(1)w'_{(4)},
w_{(3)})(w_{(1)}\otimes w_{(2)})
+\beta_{n+1; i,j}(w'_{(4)},
L(-1)w_{(3)})(w_{(1)}\otimes w_{(2)})\Big)}
\nno\\
&&=-(-n-1)\beta_{n; i,j}(w'_{(4)},
w_{(3)})(w_{(1)}\otimes w_{(2)}) 
 -(i+1)\beta_{n; i+1,j}(w'_{(4)},
w_{(3)})(w_{(1)}\otimes w_{(2)}).\nn
\end{eqnarray}
{}From (\ref{thm-11.1-10}) and (\ref{thm-11.1-12}), we obtain
\begin{eqnarray}\label{thm-11.1-12.1}
\lefteqn{((L'_{P(z_{1}-z_{2})}(0)-\wt w'_{(4)}-n-1
+\wt w_{(3)})\beta_{n; i,j}(w'_{(4)}, 
w_{(3)}))(w_{(1)}
\otimes w_{(2)})}\nno\\
&& \quad\quad\quad\quad=(\beta_{n; i,j}((L'(0)-\wt w'_{(4)})w'_{(4)}, w_{(3)}))(w_{(1)}
\otimes w_{(2)})\quad\quad\quad\quad\nno\\
&&\quad\quad\quad\quad\quad -(\beta_{n; i,j}(w'_{(4)}, (L(0)-\wt w_{(3)})w_{(3)}))(w_{(1)}
\otimes w_{(2)})\quad\quad\quad\quad\nno\\
&&\quad\quad\quad\quad\quad -(i+1)\beta_{n; i+1,j}(w'_{(4)}, (L(0)-\wt
w_{(3)})w_{(3)}))(w_{(1)}
\otimes w_{(2)}).\quad\quad\quad\quad
\end{eqnarray}
{}From (\ref{thm-11.1-12.1}), we see that there exists $N\in \Z_{+}$
independent of $w_{(1)}\in W_{(1)}$ and $w_{(2)}\in W_{(2)}$ such that
for $n\in \R$ and $i,j\in\N$,
\begin{equation}\label{thm-11.1-13}
(L'_{P(z_{1}-z_{2})}(0)-\wt w'_{(4)}+\wt w_{(3)}
-n-1)^{N}\beta_{n; i,j}(w'_{(4)}, 
w_{(3)})=0.
\end{equation}
Thus we obtain the following conclusion: For $v_{1}\in
V_{(m_{1})}^{(\alpha_{1})}, \dots, v_{r}\in
V_{(m_{r})}^{(\alpha_{r})}$, $w_{(3)}\in
(W_{3})_{[n_{3}]}^{(\beta_{3})}$ and $w'_{(4)}\in
(W_{4}')_{[n_{4}]}^{(\beta_{4})}$, the element
\[
(v^{o}_{1})_{m_{1}}\cdots 
(v^{o}_{r})_{m_{r}}\beta_{n; i,j}(w'_{(4)}, 
w_{(3)})
\]
is homogeneous of (generalized) weight
\begin{equation}\label{thm-11.1-14}
-(\mbox{\rm wt}\ v_{1}-m_{1}-1)-\cdots -(\mbox{\rm wt}\ v_{r}-m_{r}-1)
+\mbox{\rm wt}\ w'_{(4)}+n+1-\mbox{\rm wt}\ w_{(3)}
\end{equation}
and of $\tilde{A}$-degree
\[
\alpha_{1}+\cdots +\alpha_{r}+\beta_{3}+\beta_{4},
\]
and we have
\begin{eqnarray}\label{thm-11.1-15}
\lefteqn{(v^{o}_{1})_{m_{1}}\cdots (v^{o}_{r})_{m_{r}}\beta_{n; i,j}(w'_{(4)}, 
w_{(3)})=0}\nno\\
&&\quad\quad \mbox{\rm when}\;\; \wt 
(v^{o}_{1})_{m_{1}}\cdots (v^{o}_{r})_{m_{r}}\beta_{n; i,j}(w'_{(4)}, 
w_{(3)})\le N_{\alpha_{1}+\cdots +\alpha_{r}+\beta_{3}+\beta_{4}}.
\end{eqnarray}
In particular, the (generalized) weight of 
$\beta_{n; i,j}(w'_{(4)}, 
w_{(3)})$ is $\mbox{\rm wt}\ w'_{(4)}+n+1-\mbox{\rm wt}\ w_{(3)}$.

For $n\in \R$, let 
$\lambda_{n}^{(2)}(w'_{(4)}, 
w_{(3)})\in (W_{1}\otimes W_{2})^{*}$ be defined by 
\begin{eqnarray}\label{thm-11.1-16}
\lefteqn{(\lambda_{n}^{(2)}(w'_{(4)}, 
w_{(3)}))(w_{(1)}\otimes w_{(2)})}\nno\\
&&=\sum_{i,j\in\N}
(\beta_{n-\swt w'_{(4)}-1+\swt w_{(3)}; i,j}(w'_{(4)}, 
w_{(3)}))(w_{(1)}\otimes w_{(2)})\cdot\nno\\
&&\quad\quad\quad\quad\quad\cdot 
e^{(n-\swt w'_{(4)}+\swt w_{(3)})
\log (z_{1}-z_{2})}e^{(-n+\swt w'_{(4)}-\swt w_{(3)})\log z_{2}}
(\log z_{2})^{i}(\log (z_{1}-z_{2}))^{j}\nn
\end{eqnarray}
for homogeneous $w'_{(4)}\in W'_{4}$ and $w_{(3)}\in W_{3}$. 
Then 
\[
\sum_{n\in \R}(\lambda_{n}^{(2)}(w'_{(4)}, 
w_{(3)}))(w_{(1)}\otimes w_{(2)})
\]
is  absolutely convergent  to 
$\mu^{(2)}_{(I_{1}\circ (1_{W_{2}}\otimes I_{2}))'(w'_{(4)}), w_{(3)}}
(w_{(1)}\otimes w_{(2)})$ for $w_{(1)}\in W_{1}$ and $w_{(2)}\in W_{2}$.
Since  (\ref{thm-11.1-13}) holds for $n\in \R$ and $i,j\in\N$, 
\[
(L'_{P(z_{1}-z_{2})}(0)-n)^{N}\lambda_{n}^{(2)}(w'_{(4)}, 
w_{(3)})=0
\]
for $n\in \R$. 

Moreover, by (\ref{thm-11.1-16}) and (\ref{thm-11.1-13}),
\begin{eqnarray}\label{thm-11.1-17}
\lefteqn{\sum_{n\in \R}(e^{z'L'_{P(z_{1}-z_{2})}(0)}\lambda_{n}^{(2)}(w'_{(4)}, 
w_{(3)}))(w_{(1)}\otimes w_{(2)})}\nn
&&=\sum_{n\in \R}\sum_{i,j\in\N}
(e^{z'L'_{P(z_{1}-z_{2})}(0)}\beta_{n-\swt w'_{(4)}-1+\swt w_{(3)}; i,j}(w'_{(4)}, 
w_{(3)}))(w_{(1)}\otimes w_{(2)})\cdot\nno\\
&&\quad\quad\quad\quad\quad\cdot 
e^{(n-\swt w'_{(4)}+\swt w_{(3)})
\log (z_{1}-z_{2})}e^{(-n+\swt w'_{(4)}-\swt w_{(3)})\log z_{2}}
(\log z_{2})^{i}(\log (z_{1}-z_{2}))^{j}\nn
&&=\sum_{n\in \R}\sum_{i,j\in\N}\sum_{l=0}^{N-1}\frac{1}{l!}
((L'_{P(z_{1}-z_{2}}(0)-n)^{l}\beta_{n-\swt w'_{(4)}-1+\swt w_{(3)}; i,j}(w'_{(4)}, 
w_{(3)}))(w_{(1)}\otimes w_{(2)})\cdot\nno\\
&&\quad \cdot 
e^{nz'}(z')^{l}e^{(n-\swt w'_{(4)}+\swt w_{(3)})
\log (z_{1}-z_{2})}e^{(-n+\swt w'_{(4)}-\swt w_{(3)})\log z_{2}}
(\log z_{2})^{i}(\log (z_{1}-z_{2}))^{j}.\nn
\end{eqnarray}
By (\ref{thm-11.1-12.1}), for $n\in \R$, $i,j\in\N$
and $l=0, \dots, N-1$, 
\[
((L'_{P(z_{1}-z_{2}}(0)-n)^{l}\beta_{n-\swt w'_{(4)}-1+\swt w_{(3)}; i,j}(w'_{(4)}, 
w_{(3)}))(w_{(1)}\otimes w_{(2)})
\]
is a linear combination of 
\[
(\beta_{n-\swt w'_{(4)}-1+\swt w_{(3)}; i',j}((L'(0)-\wt w'_{(4)})^{p}w'_{(4)}, 
(L(0)-\wt w_{(3)})^{q}w_{(3)}))(w_{(1)}\otimes w_{(2)})
\]
for $i',p, q\in \N$ such that $p, q\le N-1$.
For $p, q\in \N$,
\begin{eqnarray*}
\lefteqn{\sum_{i,j\in\N}\sum_{n\in \R}(\beta_{n; i,j}((L'(0)-\wt w'_{(4)})^{p}w'_{(4)}, 
(L(0)-\wt w_{(3)})^{q}w_{(3)}))(w_{(1)}\otimes w_{(2)})\cdot}\nn
&&\quad\quad\quad\quad \cdot e^{(n+1)
\log (z_{1}-z_{2})}e^{(-n-1)((\log z_{2})-z')}
((\log z_{2})-z')^{i}(\log (z_{1}-z_{2}))^{j}\nn
&&=\sum_{i,j\in\N}\sum_{n\in \R}a_{n; i,j}((L'(0)-\wt w'_{(4)})^{p}w'_{(4)}, 
w_{(1)}, w_{(2)}, (L(0)-\wt w_{(3)})^{q}w_{(3)})\cdot\nn
&&\quad\quad\quad\quad \cdot e^{(\Delta+n+1)
\log (z_{1}-z_{2})}e^{(-n-1)((\log z_{2})-z')}
((\log z_{2})-z')^{i}(\log (z_{1}-z_{2}))^{j}
\end{eqnarray*}
is absolutely convergent when $z'$ is in a sufficiently small neighborhood
of $z'=0$ such that $|e^{z'}z_{2}|>|z_{0}|>0$. In particular, 
for $i,j,l\in\N$ and $p, q\in \N$ such that $p, q\le N-1$,
\begin{eqnarray*}
\lefteqn{\sum_{n\in \R}\beta_{n-\swt w'_{(4)}-1+\swt w_{(3)};
i,j}((L'(0)-\wt w'_{(4)})^{p}w'_{(4)},
(L(0)-\wt w_{(3)}^{q}w_{(3)}))(w_{(1)}\otimes w_{(2)})\cdot }\nn
&& \cdot e^{nz'}(z')^{l}e^{(n-\swt w'_{(4)}+\swt w_{(3)})
\log (z_{1}-z_{2})}e^{(-n+\swt w'_{(4)}-\swt w_{(3)})\log z_{2}}
(\log z_{2})^{i}(\log (z_{1}-z_{2}))^{j}
\end{eqnarray*}
is absolutely convergent in the same neighborhood of $z'=0$.  Thus the
right-hand side of (\ref{thm-11.1-17}) is absolutely convergent.
Hence the left-hand side of (\ref{thm-11.1-17}) is absolutely
convergent, completing the proof that $(I_{1}\circ (1_{W_{2}}\otimes
F_{2}))'(w'_{(4)})$ satisfies Part (a) of the
$P^{(2)}(z_{1}-z_{2})$-local grading restriction condition.

Recall the space $W_{\lambda}$ for $\lambda\in (W_{1}\otimes
W_{2})^{*}$ in the $P(z)$-local grading restriction condition in
Section 5 and the space $W^{(2)}_{\lambda, w_{(3)}}$ for $\lambda\in
(W_{1}\otimes W_{2}\otimes W_{3})^{*}$ and $w_{(3)}\in W_{3}$ in the
$P^{(2)}(z)$-local grading restriction condition in Section 9.  For
fixed $n\in \mathbb{R}$, $w_{(3)}\in (W_{3})_{[n_{3}]}^{(\beta_{3})}$
and $w'_{(4)}\in (W'_{4})_{[n_{4}]}^{(\beta_{4})}$,
(\ref{thm-11.1-14}) and (\ref{thm-11.1-15}) show that for $\beta\in
\tilde{A}$, the homogeneous subspace $(W_{\lambda_{n}^{(2)}(w'_{(4)},
w_{(3)})})_{[l]}^{(\beta)}$ of $W_{\lambda_{n}^{(2)}(w'_{(4)},
w_{(3)})}$ is $0$ when $l\le N_{\beta}$.  In particular, each
$\lambda_{n}^{(2)}(w'_{(4)}, w_{(3)})$, $n\in \R$, satisfies the
$P(z_{1}-z_{2})$-lower truncation condition. By Theorem \ref{9.7-1},
$W_{\lambda_{n}^{(2)}(w'_{(4)}, w_{(3)})}$ is a doubly-graded
generalized $V$-module (a doubly-graded $V$-module when ${\cal C}$ is
in ${\cal M}_{sg}$).  Since for $\beta\in \tilde{A}$,
$(W_{\lambda_{n}^{(2)}(w'_{(4)}, w_{(3)})})_{[l]} ^{(\beta)}=0$ when
$l\le N_{\beta}$, it is in fact lower bounded.  Thus
$W_{\lambda_{n}^{(2)}(w'_{(4)}, w_{(3)})}$, generated by
$\lambda_{n}^{(2)}(w'_{(4)}, w_{(3)})$, is a finitely-generated lower
bounded doubly-graded generalized $V$-module (or $V$-module), and so
by assumption it is in fact an object of $\mathcal{C}$. Thus
$W^{(2)}_{(I_{1}\circ (1_{W_{2}}\otimes I_{2}))' (w'_{(4)}),
w_{(3)}}$, as a sum of these objects of $\mathcal{C}$, lies in
$W_{1}\hboxtr_{P(z_{1}-z_{2})}W_{2}$, which, by assumption, is an
object of $\mathcal{C}$.  Since $\mathcal{C}$ is a subcategory of
$\mathcal{GM}_{sg}$, $W^{(2)}_{(I_{1}\circ (1_{W_{2}}\otimes I_{2}))'
(w'_{(4)}), w_{(3)}}$, as a generalized $V$-module of
$W_{1}\hboxtr_{P(z_{1}-z_{2})}W_{2}$, satisfies the two
grading-restriction conditions, proving that the element $I_{1}\circ
(1_{W_{2}}\otimes I_{2}))'(w'_{(4)})$ of $(W_{1}\otimes W_{2}\otimes
W_{3})^{*}$ satisfies the $P^{(2)}(z_{1}-z_{2})$-local
grading-restriction condition.  \epfv

\subsection{Differential equations}

In this section, we assume for simplicity that $A$ and $\tilde{A}$ are
trivial.  We would like to emphasize that this assumption is not
essential since all the results can be generalized to the case that
$A$ and $\tilde{A}$ are not trivial.  To avoid spending too many pages
to straighforwardly generalize many definitions and results to the
general case (for example, the definition of $C_{1}$-cofiniteness),
and to allow the use of most of the arguments in \cite{diff-eqn}, for
which both $A$ and $\tilde{A}$ are trivial, we choose to simply assume
that $A$ and $\tilde{A}$ are trivial.

We use differential equations to prove the convergence and extension
property, given in the preceding section, when the objects of
$\mathcal{C}$ satisfy natural conditions.  The results and the proofs
here are in fact the same as those in \cite{diff-eqn}, except that in
this section, we consider objects of $\mathcal{C}$, not just ordinary
$V$-modules, and we consider logarithmic intertwining operators, not
just ordinary intertwining operators. So in the proofs of the results
in this section, we shall indicate only how the proofs here differ
{}from the corresponding ones in \cite{diff-eqn} and refer the reader
to \cite{diff-eqn} for more details.

Let $V$ be a M\"{o}bius or conformal vertex algebra with $A$ the
trivial group and let
\[
V_{+}=\coprod_{n>0}V_{(n)}.
\]
For a generalized $V$-module $W$, set
\[
C_{1}(W)={\rm span}\{u_{-1}w\;|\; u\in V_{+},\;
w\in W\}.  
\]

\begin{defi}
{\rm If $W/C_{1}(W)$ is finite dimensional, we say that $W$ is {\it
$C_{1}$-cofinite} or satisfies the {\it $C_{1}$-cofiniteness
condition}.  If for any $N\in {\mathbb R}$,
$\coprod_{n<N}W_{[n]}$ is finite dimensional, we say that $W$
is {\it quasi-finite dimensional} or satisfies the {\it
quasi-finite-dimensionality condition}.}
\end{defi}

We have:

\begin{theo}\label{sys}
Let $W_{i}$ for $i=0, \dots, n+1$ be generalized $V$-modules
satisfying the $C_{1}$-cofiniteness condition and the
quasi-finite-dimensionality condition.  Then for any $w'_{(0)}\in
W'_{0}$, $w_{(1)}\in W_{1}, \dots, w_{(n+1)}\in W_{n+1}$, there exist
\[
a_{k, \;l}(z_{1}, \dots, z_{n})\in {\mathbb C}[z_{1}^{\pm 1}, \dots,
z_{n}^{\pm 1}, (z_{1}-z_{2})^{-1}, (z_{1}-z_{3})^{-1}, \dots,
(z_{n-1}-z_{n})^{-1}],
\]
for $k=1, \dots, m$ and $l=1, \dots, n,$ such that the following
holds: For any generalized $V$-modules $\widetilde{W}_{1}, \dots,
\widetilde{W}_{n-1}$, and any logarithmic intertwining operators
\[
{\cal Y}_{1}, {\cal Y}_{2}, \dots, {\cal
Y}_{n-1}, {\cal Y}_{n}
\]
of types 
\[
{W_{0}\choose
W_{1}\widetilde{W}_{1}}, {\widetilde{W}_{1}\choose
W_{2}\widetilde{W}_{2}}, \dots, {\widetilde{W}_{n-2}\choose
W_{n-1}\widetilde{W}_{n-1}}, {\widetilde{W}_{n-1}\choose
W_{n}W_{n+1}},
\]
respectively, the series
\begin{equation}
\langle w'_{(0)}, {\cal Y}_{1}(w_{(1)}, z_{1})\cdots {\cal Y}_{n}(w_{(n)},
z_{n})w_{(n+1)}\rangle
\end{equation}
satisfies the system of differential equations
\[
\frac{\partial^{m}\varphi}{\partial z_{l}^{m}}+ \sum_{k=1}^{m}
\iota_{|z_{1}|>\cdots >|z_{n}|>0}(a_{k,
\;l}(z_{1}, \dots, z_{n})) \frac{\partial^{m-k}\varphi}{\partial
z_{l}^{m-k}}=0,\;\;\;l=1, \dots, n
\]
in the region $|z_{1}|>\cdots >|z_{n}|>0$, where 
\[
\iota_{|z_{1}|>\cdots >|z_{n}|>0}(a_{k,
\;l}(z_{1}, \dots, z_{n}))
\]
for $k=1, \dots, m$ and $l=1, \dots, n$ are the (unique) Laurent
expansions of $a_{k, \;l}(z_{1}, \dots, z_{n})$ in the region
$|z_{1}|>\cdots >|z_{n}|>0$.  Moreover, for any set of possible
singular points of the system
\begin{equation}\label{sys-eqns}
\frac{\partial^{m}\varphi}{\partial z_{l}^{m}}+ \sum_{k=1}^{m} a_{k,
\;l}(z_{1}, \dots, z_{n})\frac{\partial^{m-k}\varphi}{\partial
z_{l}^{m-k}}=0,\;\;\;l=1, \dots, n
\end{equation}
such that either $z_i = 0$ or $z_i = \infty$ for some $i$ or $z_i =
z_j$ for some $i \ne j$, the $a_{k, \;l}(z_{1}, \dots, z_{n})$ can be
chosen for $k=1, \dots, m$ and $l=1, \dots, n$ so that these singular
points are regular.
\end{theo}
\pf Proposition 1.1, Corollary 1.2 and Corollary 1.3 in
\cite{diff-eqn} and their proofs still hold when all $V$-modules are
replaced by strongly-graded generalized $V$-modules. Theorem 1.4 in
\cite{diff-eqn} also still holds, except that in its proof all
$V$-modules are replaced by strongly-graded generalized $V$-modules
and the spaces
\begin{eqnarray*}
&z_{1}^{\Delta}
\C(\{z_{2}/z_{1}\})[z_{1}^{\pm 1}, z_{2}^{\pm 1}],&\\
&z_{2}^{\Delta}
\mathbb{C}(\{(z_{1}-z_{2})/z_{1}\})[z_{2}^{\pm 1},
(z_{1}-z_{2})^{\pm 1}],&\\
&z_{2}^{\Delta}
\C(\{z_{1}/z_{2}\})[z_{1}^{\pm 1}, z_{2}^{\pm 1}]&
\end{eqnarray*}
are replaced by 
\begin{eqnarray*}
&z_{1}^{\Delta}
\C(\{z_{2}/z_{1}\})[z_{1}^{\pm 1}, z_{2}^{\pm 1}]
[\log z_{1}, \log z_{2}],&\\
&z_{2}^{\Delta}
\mathbb{C}(\{(z_{1}-z_{2})/z_{1}\})[z_{2}^{\pm 1},
(z_{1}-z_{2})^{\pm 1}][\log z_{2}, \log (z_{1}-z_{2})],&\\
&z_{2}^{\Delta}
\C(\{z_{1}/z_{2}\})[z_{1}^{\pm 1}, z_{2}^{\pm 1}]
[\log z_{1}, \log z_{2}],&
\end{eqnarray*}
respectively.  Theorem 1.6 and its proof also still hold in our
setting here.  Thus the first conclusion holds.

Proposition 2.1, Lemma 2.2, Theorem 2.3, Remark 2.4 and Theorem 2.5 in
\cite{diff-eqn} and their proofs still hold here.  Thus the second
conclusion holds.  \epfv

\begin{rema}\label{simple-sing}
{\rm In \cite{diff-eqn} and in Theorem \ref{sys} above, we are using
the definition of the notion of system of differential equations with
regular singular points (or regular singularities) given in Appendix B
of the book \cite{Kn}. This definition defines such a system in terms
of the form of the expansions of the solutions at the singular
points. What we need in the present section is precisely this form of
the expansions of the solutions of such a system.  However, we warn
the reader that in some books and papers, a system of differential
equations with regular singular points is instead defined in terms of
the coefficients of the system.  In fact, there are systems with
regular singular points that do not satisfy the definition in terms of
the coefficients of the system.  See for example \cite{Lut} for a
study of the problem of analyzing under what conditions the singular
points of a system of differential equations that do not satisfy the
regularity condition in terms of the coefficients are still regular.
In Section B.5 of \cite{Kn}, systems with simple singular points
(called simple sigularities there and defined in terms of the
coefficients) are discussed and it is proved there that simple
singular points are indeed regular singular points. From Theorem B.16
in Section B.5 of \cite{Kn}, together with (2.2) and (2.3) in
\cite{diff-eqn} and their obvious extensions to the case of more than
two intertwining operators, it is clear that for any set of possible
of singular points of (\ref{sys-eqns}) of any of the indicated types,
the $a_{k, \;l}(z_{1}, \dots, z_{n})$ can be chosen for $k=1, \dots,
m$ and $l=1, \dots, n$ so that these singular points are simple
singular points and thus are regular singular points.}
\end{rema}

We now have:

\begin{theo}\label{C_1pp}
Suppose that all generalized $V$-modules in ${\cal C}$ satisfy the
$C_{1}$-cofiniteness condition and the quasi-finite-dimensionality
condition.  Then:

\begin{enumerate}

\item The convergence and extension properties for products and
iterates hold in ${\cal C}$. If $\mathcal{C}$ is in $\mathcal{M}_{sg}$
and if every object of $\mathcal{C}$ is a direct sum of irreducible
objects of $\mathcal{C}$ and there are only finitely many irreducible
objects of $\mathcal{C}$ (up to equivalence), then the convergence and
extension properties without logarithms for products and iterates hold
in ${\cal C}$.

\item For any $n\in{\mathbb Z}_+$, any objects $W_1, \dots, W_{n+1}$
and $\widetilde{W}_1, \dots, \widetilde{W}_{n-1}$ of $\mathcal{C}$,
any logarithmic intertwining operators
\[
{\cal Y}_{1}, {\cal Y}_{2}, \dots, {\cal
Y}_{n-1}, {\cal Y}_{n}
\]
of types 
\[
{W_{0}\choose
W_{1}\widetilde{W}_{1}}, {\widetilde{W}_{1}\choose
W_{2}\widetilde{W}_{2}}, \dots, {\widetilde{W}_{n-2}\choose
W_{n-1}\widetilde{W}_{n-1}}, {\widetilde{W}_{n-1}\choose
W_{n}W_{n+1}},
\]
respectively, and any $w_{(0)}'\in W_{0}'$, $w_{(1)}\in W_{1}, \dots,
w_{(n+1)}\in W_{n+1}$, the series
\begin{equation} 
\langle
w_{(0)}', {\cal Y}_{1}(w_{(1)}, z_{1})\cdots {\cal Y}_{n}(w_{(n)},
z_{n})w_{(n+1)}\rangle
\end{equation}
is absolutely convergent in the region $|z_{1}|>\cdots> |z_{n}|>0$ and
its sum can be analytically extended to a multivalued analytic
function on the region given by $z_{i}\ne 0$, $i=1, \dots, n$,
$z_{i}\ne z_{j}$, $i\ne j$, such that for any set of possible singular
points with either $z_{i}=0$, $z_{i}=\infty$ or $z_{i}= z_{j}$ for
$i\ne j$, this multivalued analytic function can be expanded near the
singularity as a series having the same form as the expansion near the
singular points of a solution of a system of differential equations
with regular singular points (as defined in Appendix B of \cite{Kn};
recall Remark \ref{simple-sing}).

\end{enumerate}

\end{theo}
\pf The first statement in the first part, except for the existence of
the upper bound $K$ for $i_{k}$, and the statement in the second part
of the theorem follow directly {}from Theorem \ref{sys} and the theory
of differential equations with regular singular points.  To prove the
existence of $K$, note that in \cite{diff-eqn}, $T/J$ is a finitely
generated $R$-module and the map $\phi_{\Y_{1}, \Y_{2}}$ with domain
$T$ induces a map $\bar{\phi}_{\Y_{1}, \Y_{2}}$ with domain $T/J$ such
that in particular the image of $\phi_{\Y_{1}, \Y_{2}}$ is equal to
the image of $\bar{\phi}_{\Y_{1}, \Y_{2}}$.  Since $T/J$ is a
finitely-generated $\C[z_{1}^{\pm 1}, z_{2}^{\pm 1},
(z_{1}-z_{2})^{-1}]$-module, its image under $\bar{\phi}_{\Y_{1},
\Y_{2}}$ is also a finitely-generated $\C[z_{1}^{\pm 1}, z_{2}^{\pm
1}, (z_{1}-z_{2})^{-1}]$-module.  Any set of generators of this image
has been proved to have analytic extensions in the region
$|z_{2}|>|z_{1}-z_{2}|>0$ of the form (\ref{c-e-p-3}), and for
finitely many such generators, we have an upper bound $K$ for the
powers of $\log z_{2}$ in their analytic extensions in this region.
Thus the powers of $\log z_{2}$ in the analytic extension of any
element of the image of the map $\phi_{\Y_{1}, \Y_{2}}$ in this
region, as a linear combination of the analytic extensions of these
generators with coefficients in $\C[z_{1}^{\pm 1}, z_{2}^{\pm 1},
(z_{1}-z_{2})^{-1}]$, are also less than $K$.

The second statement in the first part was in fact proved in
\cite{diff-eqn}.  \epfv

\begin{rema}
{\rm Note that the first statement in the first part of Theorem
\ref{C_1pp} follows immediately from Theorem \ref{sys}.
To prove the
second statement in the first part of Theorem \ref{C_1pp}, we start
with the product of two intertwining operators without logarithms in
the region $|z_{1}|>|z_{2}|>0$, and we have to prove that its analytic
extension in the region $|z_{2}|>|z_{1}-z_{2}|>0$ has no terms
involving logarithms. This is the main hard part of the proof of
Theorem 3.5 in \cite{diff-eqn}.}
\end{rema}


\bigskip

\noindent {\small \sc Department of Mathematics, Rutgers University,
Piscataway, NJ 08854 (permanent address)}

\noindent {\it and}

\noindent {\small \sc Beijing International Center for Mathematical Research,
Peking University, Beijing, China}

\noindent {\em E-mail address}: yzhuang@math.rutgers.edu

\vspace{1em}

\noindent {\small \sc Department of Mathematics, Rutgers University,
Piscataway, NJ 08854}

\noindent {\em E-mail address}: lepowsky@math.rutgers.edu

\vspace{1em}

\noindent {\small \sc Department of Mathematics, Rutgers University,
Piscataway, NJ 08854}

\noindent {\em E-mail address}: linzhang@math.rutgers.edu

\end{document}